\documentclass{ifacconf}

\usepackage{graphicx}      
\usepackage{natbib}        
\usepackage{amsmath}
\usepackage{amssymb}
\usepackage{xcolor}
 \usepackage{epsfig}
\usepackage[normalsize,it,hang]{subfigure}
\usepackage{lipsum}

\begin{document}

\begin{frontmatter}

{\tt MATHMOD 2025 \\ 11th Vienna International Conference on \\Mathematical Modelling \\
19-21 February 2025, Vienna, Austria}

\title{Synchronization of Kuramoto oscillators \\ via HEOL, and a  discussion on AI}

\author[First]{Emmanuel Delaleau} 
\author[Second,fourth]{C\'{e}dric Join}
\author[Third,fourth]{Michel Fliess}

\address[First]{ENI Brest, IRDL (CNRS, UMR 6027), 29200 Brest, France\\
  (e-mail: emmanuel.delaleau@enib.fr)}

\address[Second]{CRAN (CNRS, UMR 7039), Universit\'{e} de Lorraine,
  BP 239,54506 Vand{\oe}uvre-l\`{e}s-Nancy, France (e-mail:
  cedric.join@univ-lorraine.fr)}

\address[Third]{LIX (CNRS, UMR 7161), \'Ecole polytechnique, 91128
  Palaiseau, France (e-mail:
  michel.fliess@polytechnique.edu, michel.fliess@swissknife.tech)}

\address[fourth]{AL.I.E.N., 7 rue Maurice Barr\`{e}s, 54330 V\'{e}zelise, France \\
(e-mail: \{cedric.join, michel.fliess\}@alien-sas.com)}

\begin{abstract}
Artificial neural networks and their applications in deep learning have recently made an incursion into the field of control. Deep learning techniques in control are often related to optimal control, which relies on the Pontryagin maximum principle or the Hamilton-Jacobi-Bellman equation. They imply control schemes that are tedious to implement. We show here that the new HEOL setting, resulting from the fusion of the two established approaches, namely differential flatness and model-free control, provides a solution to control problems that is more sober in terms of computational resources.
This communication is devoted to the synchronization of the popular Kuramoto's coupled oscillators, which was already considered via artificial neural networks by L. Böttcher \emph{et al.} (Nature Commun., 2022), where, contrarily to this communication, only the single control variable case is examined. 
One establishes the flatness of Kuramoto’s coupled oscillator model with multiplicative control and develops the resulting HEOL control. Unlike many examples, this system reveals singularities that are avoided by a clever generation of phase angle trajectories. The results obtained, verified in simulations, show that it is not only possible to synchronize these oscillators in finite time, and even to follow angular frequency profiles, but also to exhibit robustness concerning model mismatches. To the best of our knowledge that has never been done before. Concluding remarks advocate a viewpoint, which might be traced back to Wiener's cybernetics: control theory belongs to AI.
\end{abstract}


\end{frontmatter}

\onecolumn


\section{Introduction}

The stunning advances in the field of deep learning via artificial neural networks (ANNs) (see, e.g., \citet{cun}) explain why, contrarily to the situation some time ago (see, e.g., \cite{sutton}), the relationship between control engineering and artificial intelligence is today often investigated with ANNs (see, e.g., among a huge number of publications, Narendra and Parthasarathy (1990), \cite{suy}, \cite{sara}, \cite{dev}, \cite{ben}, \cite{cerf}, \cite{zhou}). This communication starts analyzing this situation by way of a recent paper due to \cite{Bottcher22}, where the synchronization of some oscillators introduced by \cite{Kuramoto75,Kuramoto84} is considered through ANNs. The popularity of those oscillators, the mathematical modeling of which is known, is explained by their surprising versatility in a variety of fields, all seemingly foreign to one another, ranging from pure physics and chemical reactions to smart grids and neurosciences: See, e.g., \cite{Acebron05}, \cite{chopra}, \cite{Breakspear10}, \cite{Dorfler14}, \cite{Dorfler13}, \cite{lss}, \cite{Stroglatz05}, and references therein. They bear witness to the diversity of the techniques used, including those of control.

In our approach to synchronization, we follow \cite{Bottcher22} for placing multiplicative control variables.  It permits to use of some tools that seem to have never been applied before for this type of question:
\begin{enumerate}
    \item The multivariable controlled system becomes obviously (\emph{differentially}) \emph{flat} (Fliess et al. (1995, 1999)). The phase angles are the \emph{flat outputs}. This notion, which is now quite popular in engineering (see, e.g., the books by \cite{hagen}, \cite{sira}, \cite{levine}, \cite{rudolph}), yields open-loop reference trajectories for the phase angles, which not only ensure synchronicity in finite time but also a convenient behavior for the system variables. To the best of our knowledge, such results were not achieved until now.
    \item The loop is closed via the \emph{HEOL} setting (\cite{heol}), which is inspired from \emph{model-free control} (Fliess and Join (2013, 2022)): See references there and in \cite{heol} for numerous examples of successful concrete applications. Model mismatches are therefore easily handled. To the best of our knowledge, such robustness issues seem not to have been investigated in the existing literature. 
\end{enumerate}
We also briefly examine the modeling proposed by \cite{Mao16} where the control variables are additive and not multiplicative like in \cite{Bottcher22}. It is again trivially flat.

It is well known that today's machine learning techniques are most often intimately related to techniques stemming from optimal control, like the Pontryagin maximum principle and the Hamilton-Jacobi-Bellman partial differential equation which are in general most difficult to implement in practice despite many attempts (see, e.g., \cite{mit} and \cite{jin}), especially in concrete control engineering. Our results confirm therefore \cite{vancouver}, which was about model-free control. Appropriate theoretical advances in control engineering seem to perform better at least in the continuous-time case today than deep learning via artificial neural networks. 

In a most original contribution \cite{Bottcher22} are introducing \emph{AI Pontryagin}, i.e., ANN methods to bypass the tedious calculations related to optimal control. Excellent computer experiments are depicted. A thorough comparison seems difficult to develop here in such a restricted place. Let us emphasize however that \cite{Bottcher22} is only dealing with the single control variable case, contrary to what is presented here. It is also quite clear that from the point of view of computational power AI Pontryagin is much more demanding than our HEOL setting. See, e.g., \cite{chaxel} for the feedback implementation. 

Our paper is organized as follows. Sect. \ref{flatness} is devoted to the flatness-based open-loop control of the modeling due to \cite{Bottcher22}. The closed-loop control is examined in Sect. \ref{HEOL}. Several computer simulations are examined in Sect. \ref{exp}. Additive control variables are briefly treated in Sect. \ref{add}. Some concluding remarks on future investigations and on the relationship with AI may be found in Sect. \ref{conclu}.

\section{Flatness-based open-loop control}\label{flatness}

\cite{Bottcher22} consider the following Kuramoto model of coupled oscillators:
\begin{equation}
  \label{eq:model:kuramoto}
  \dot\theta_i =
  \omega_i +
 u_i \frac{K}{N}\sum_{j=1}^N a_{i,j} \sin(\theta_j- \theta_i), \quad i = 1,\ldots,N
\end{equation}
where
$N$ is the number of oscillators, $\theta_i$, $\omega_i$ and $u_i$ are respectively the phase angle, the natural angular frequency, and the control variable associated with the $i$th oscillator, $K$ is the coupling strength,  the $a_{i,j}$'s are adjacency coefficients. It is obvious that Eq.~\eqref{eq:model:kuramoto} defines a flat system, where $\theta_i$, $i = 1, \dots, N$, is a flat output (Fliess et al. (1995, 1999)). The control variables $u_i$ are functions of the flat outputs and their derivatives:
\begin{equation}\label{flat}
   u_i =  \frac{N(\dot{\theta}_i - \omega_i)}{K\sum_{j=1}^N a_{i,j} \sin(\theta_j - \theta_i)}, \quad i = 1,\ldots,N
\end{equation}


\subsection{Reference trajectories}\label{sing}
The main difficulties in assigning reference trajectories to the flat outputs for achieving synchronization, i.e., ${\dot \theta}_1 (t) = \dots = {\dot\theta}_N (t)$, for $t \geqslant t_f$, are the following:
\begin{enumerate}
    \item The denominators in Eq. \eqref{flat} should not be equal to $0$, i.e., $\sum_{j=1}^N a_{i,j} \sin(\theta_j(t) - \theta_i(t)) \neq 0$, $\forall t$.
    \item The derivatives of the phases should be positive, i.e., ${\dot \theta}_i(t) > 0$, $\forall t$.
    \item The control variables $u_i$ should be positive, $\forall t$.
\end{enumerate}

Write
$$ \theta_i (t) = g_i(t) + f(t), \quad \forall i=1,\ldots,N$$
where $f$ is the synchronization function.

\begin{figure*}[!ht]
\centering
\subfigure[\footnotesize $\delta \theta_1$]
{\epsfig{figure=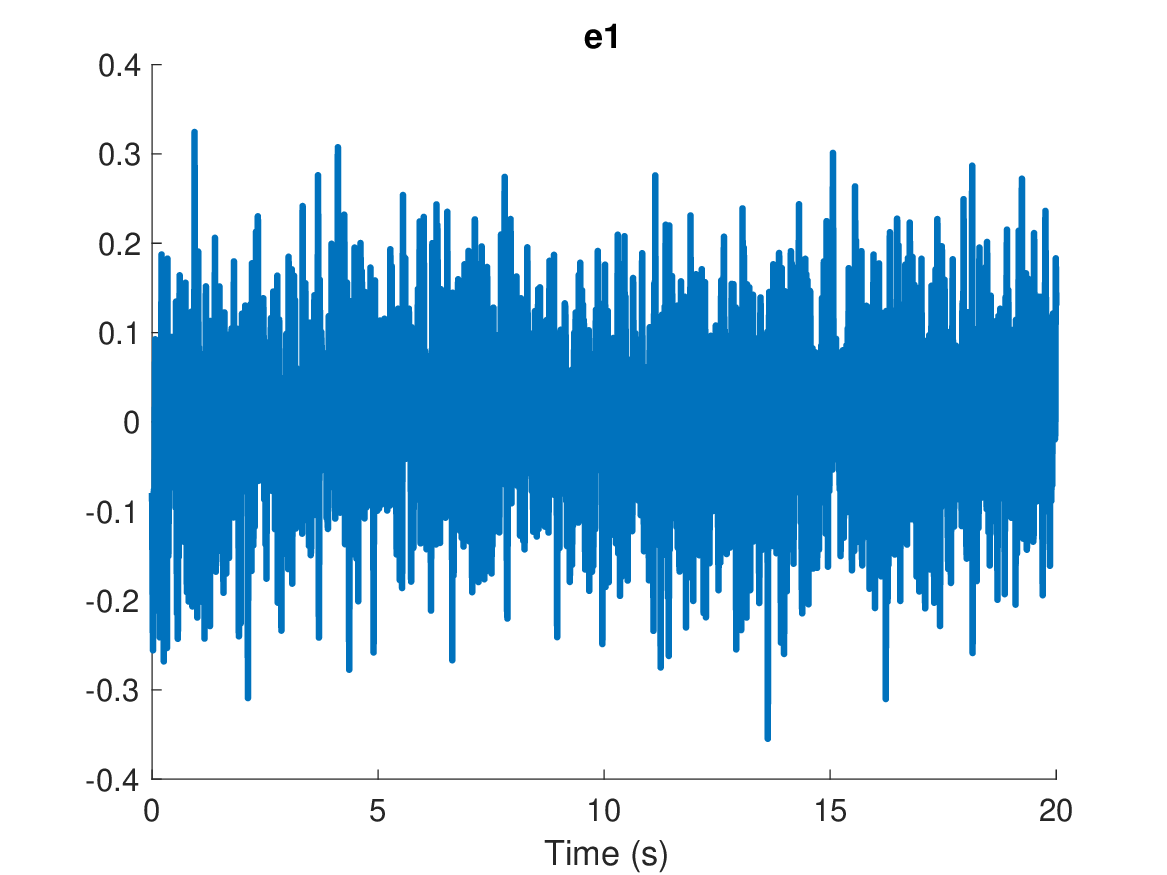,width=0.325\textwidth}}
\subfigure[\footnotesize $\delta \theta_2$]
{\epsfig{figure=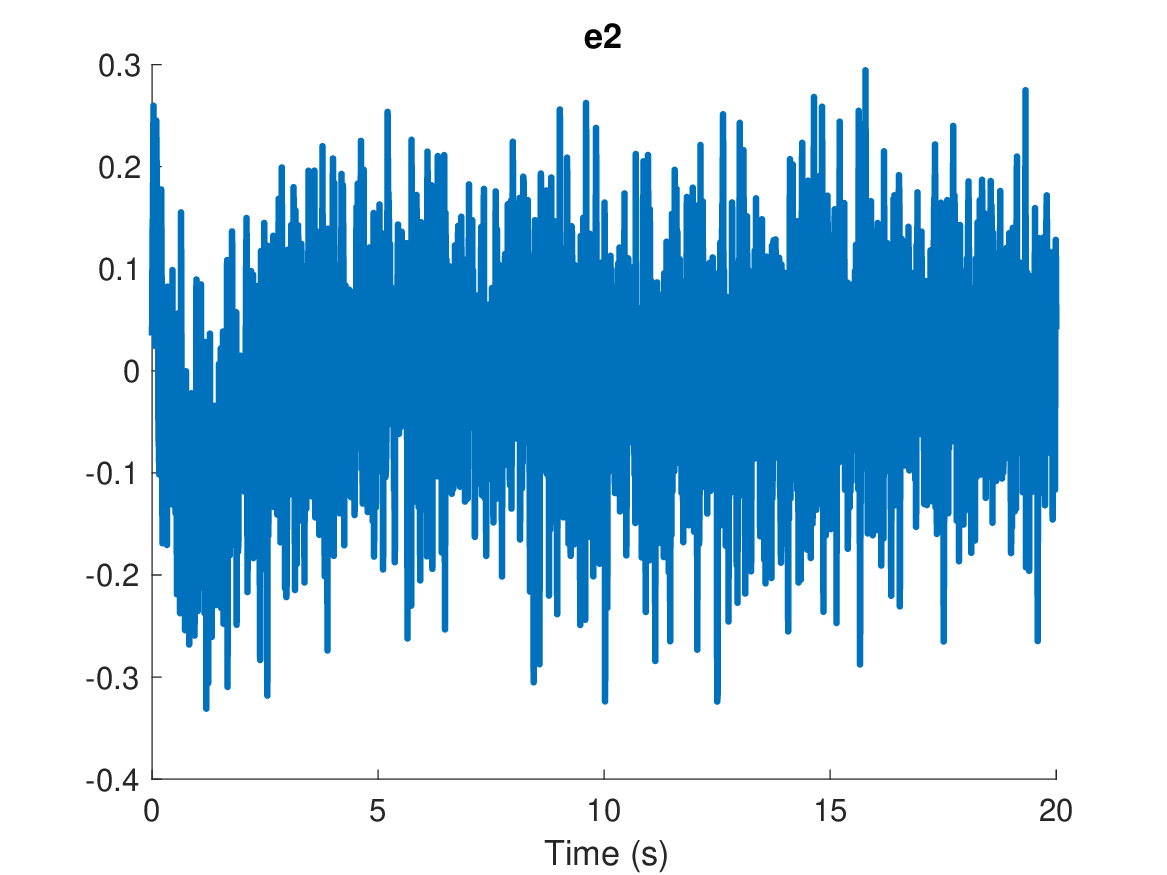,width=0.325\textwidth}}
\subfigure[\footnotesize $\delta \theta_3$]
{\epsfig{figure=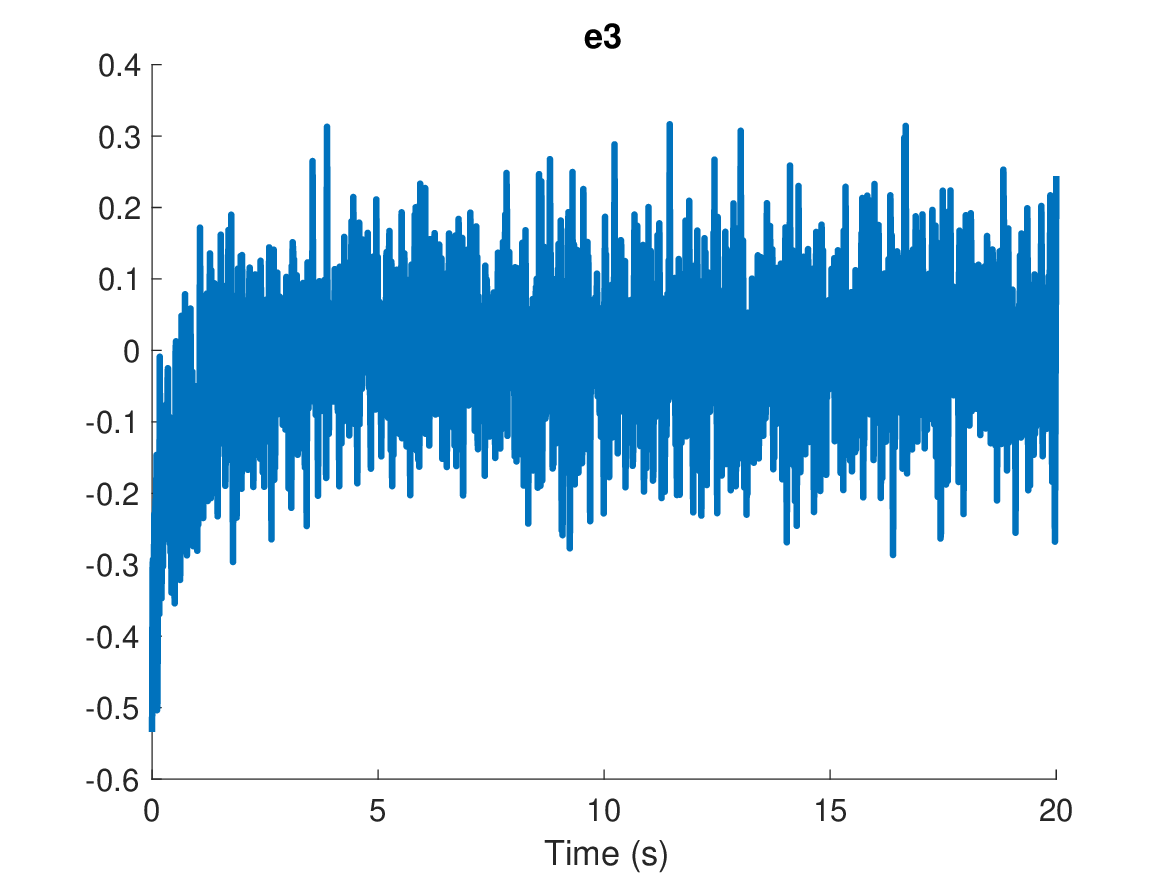,width=0.325\textwidth}}
\caption{Tracking errors}\label{Be}
\end{figure*}

The synchronization, i.e., $\dot{\theta}_i(t) \approx \dot{\theta}_j(t)$, $i\neq j$, for $t \geqslant t_f$, is equivalent to the fact that the $g_i(t)$'s are approximately constant for $t \geqslant t_f$. The following linear differential equation
\begin{equation}
    \label{eq:diff:gi}
    \tau^2 \ddot g_i +2\tau \dot g_i+g_i=c_i
\end{equation}
 easily achieves this. The solution reads:
$g_i(t) = c_i + (A + Bt)\exp\left(-\frac{t}{\tau}\right)$
where $A$ and $B$ are constants, which are deduced from the initial conditions. One can impose, e.g.,  that $\left|\frac{g_i(t_f) - c_i}{c_i}\right| \leqslant 0.001$ for a large enough $t_f$, i.e., $t_f$ is chosen to be equal to a few time constants~$\tau$.
The generation of trajectories as the output of a second-order filter has already been used in electrical drives \citep{DelalHagen02jesa}. See Sect. \ref{exp} and Figs. \ref{By} to \ref{Bey}. 

\section{Closed-loop control via HEOL}\label{HEOL}

\subsection{The homeostat}

Differentiate Eq. \eqref{eq:model:kuramoto}: 
\begin{eqnarray}\label{lin}
 d\dot\theta_i &=& d u_i \frac{K}{N}\sum_{j=1}^N a_{i,j} \sin(\theta_j- \theta_i) \nonumber\\
 &&\mbox{}+u_i \frac{K}{N}\sum_{j=1}^N a_{i,j} (d\theta_j - d\theta_j)\cos(\theta_j- \theta_i) \label{eq:model:diffkuramoto}\\
 \nonumber &&i = 1,\ldots,N
\end{eqnarray}

In the HEOL\footnote{\emph{Sun} in the Breton language.} setting \citep{heol}, Eq.~\eqref{lin} should be understood  as the \emph{homeostat}.\footnote{Terminology borrowed from \cite{ashby}.} Write $u_i^\star$ and $\theta_i^\star$ the control variables and the corresponding reference trajectories, and $\delta u_i = u_i - u_{i}^{\star}$, $\delta \theta_i = \theta_i - \theta_{i}^{\star}$, $\delta\dot\theta_{i} = \dot{\theta}_i - \dot{\theta}_{i}^{\star} = \frac{d}{dt} \delta \theta_i$,
\begin{equation}\label{eq:homeo}
    \frac{d}{dt}\delta \theta_i = F_i +\alpha_i \delta u_i
\end{equation}
where
\begin{eqnarray*}
    F_i &=& u_{i}^{\star}\frac{K}{N}
    \sum_{j=1}^N a_{i,j}
    (\delta \theta_j - \delta \theta_i)\cos(\theta_{j}^{\star} - \theta_{i}^{\star} )
    \\
    \alpha_i &=& \frac{K}{N}\sum_{j=1}^N a_{i,j} \sin(\theta_{j}^{\star}  - \theta_{i}^{\star} )
\end{eqnarray*}
Following \citet{heol}, $F_i$ stands now for the mismatches and disturbances, like in the well-known \emph{ultra-local} model of model-free control (Fliess and Join (2013, 2022)). But contrarily to the classic model-free control approach, the coefficient $\alpha_i$ of $\delta u_i$ may be time-varying.

Techniques from operational calculus yield a \emph{data-driven real-time} estimator $F_{i}^{\rm{est}}$ \citep{heol} of $F_i$:
 
\begin{eqnarray}\label{estim}
 F_{i}^{\rm{est}} &=& - \frac{6}{T^3} \int_{0}^{T} ( (T - 2 \sigma)\delta \theta_i(\sigma+t-T) 
 \\ \nonumber
 &+& \sigma (T - \sigma)\alpha_i(\sigma+t-T)\delta u_i(\sigma+t-T))d\sigma
\end{eqnarray}
where $T > 0$ is ``small.'' In the parlance of today's AI Formula \eqref{estim} might be
viewed as a peculiar type of machine learning.

The corresponding \emph{intelligent proportional controller}, or \emph{iP}, reads
\begin{equation}\label{ip}
  \delta u_i = -\frac{F_{i}^{\rm{est}} + K_{P, i} \delta \theta_i}{\alpha_i}  
\end{equation}
where $K_{P, i}$ is the gain. Combining Eq. \eqref{eq:homeo} and \eqref{ip} yields
$$\frac{d}{dt} \delta \theta_i + K_{P, i} \delta \theta_i = F_i - F_{i}^{\rm{est}} $$
Assume that the estimate $F_{i}^{\rm{est}}$ is ``good,'' i.e., $F_i - F_{i}^{\rm{est}} \approx 0$, then $\lim_{t \rightarrow + \infty} \delta \theta_i (t) \approx 0$ if $K_{P, i} > 0$. Local stability around the reference trajectory is ensured.

\subsection{Computer experiments}\label{exp}

\begin{figure*}[!ht]
\centering
\subfigure[\footnotesize $u_1$ (blue--), $u_1^\star$ (red - -)]
{\epsfig{figure=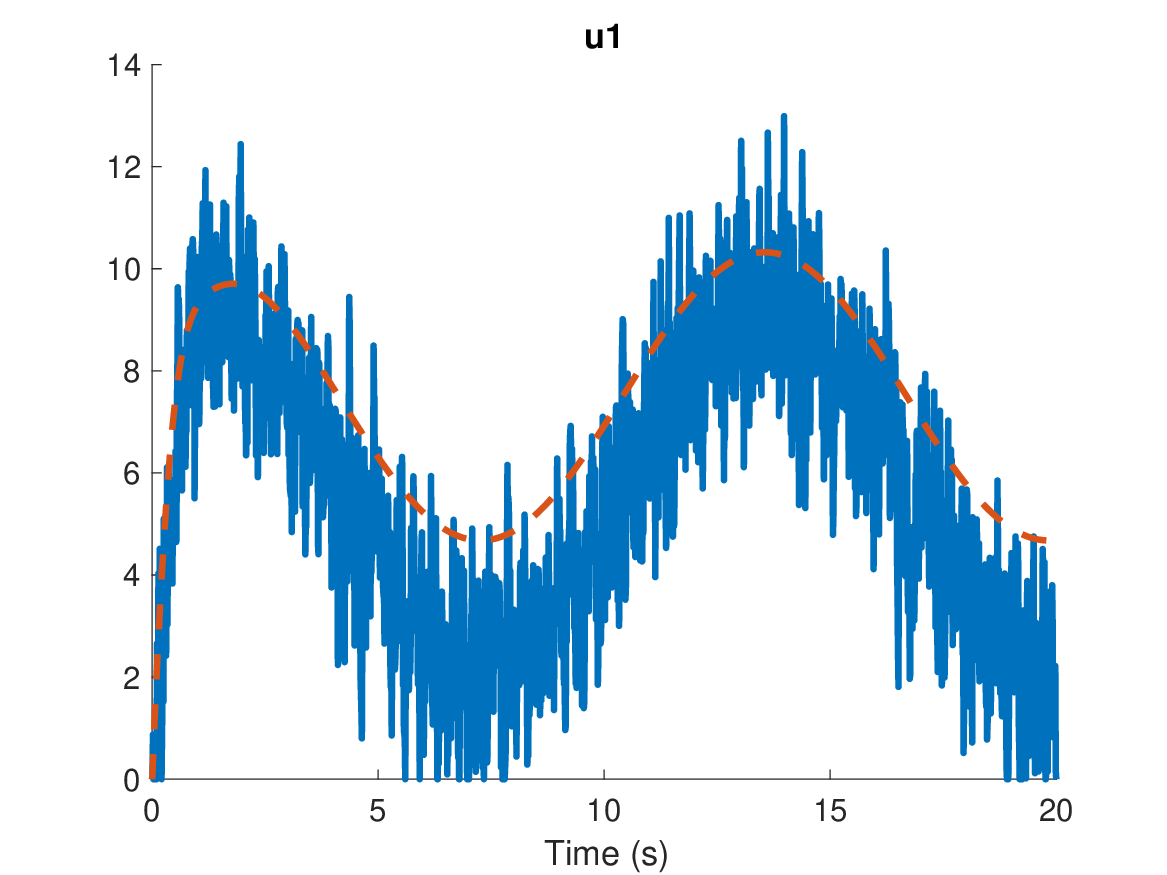,width=0.325\textwidth}}
\subfigure[\footnotesize $u_2$ (blue--), $u_2^\star$ (red - -)]
{\epsfig{figure=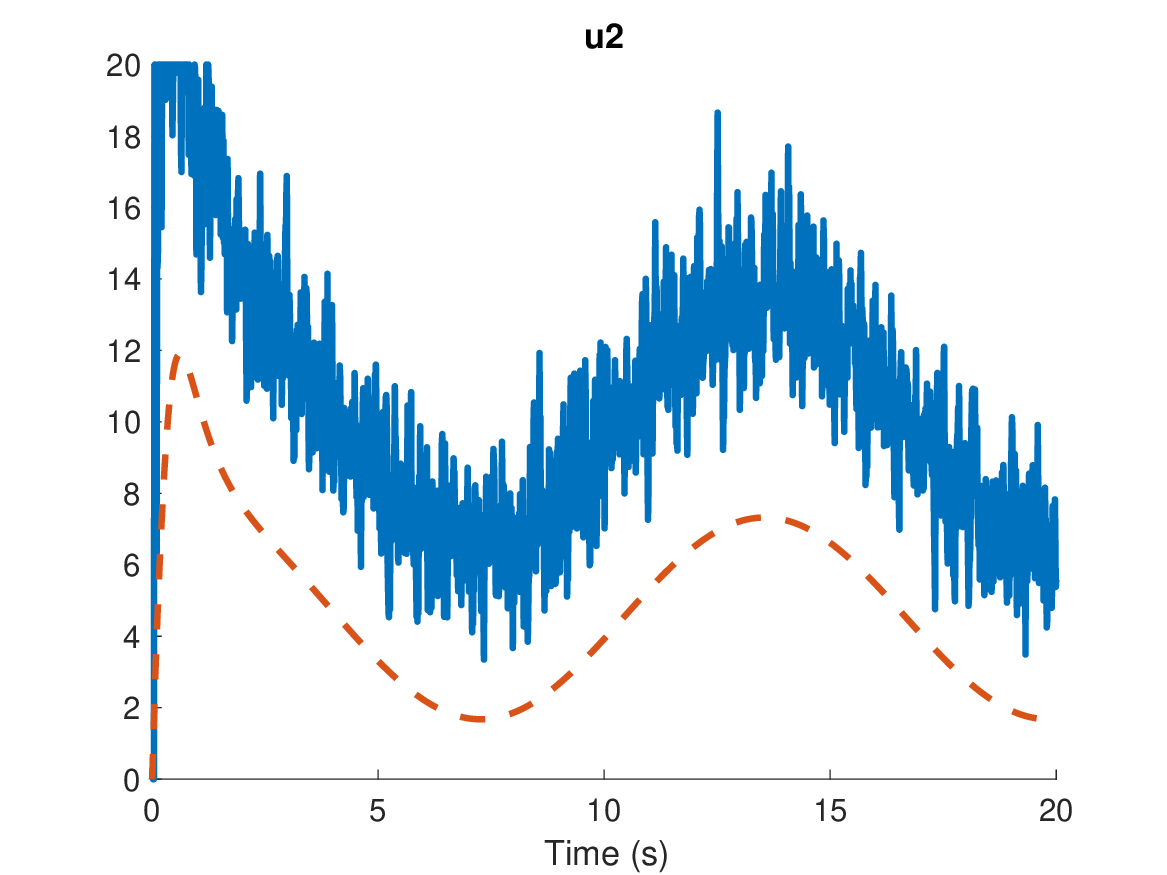,width=0.325\textwidth}}
\subfigure[\footnotesize$u_3$ (blue--), $u_3^\star$ (red - -)]
{\epsfig{figure=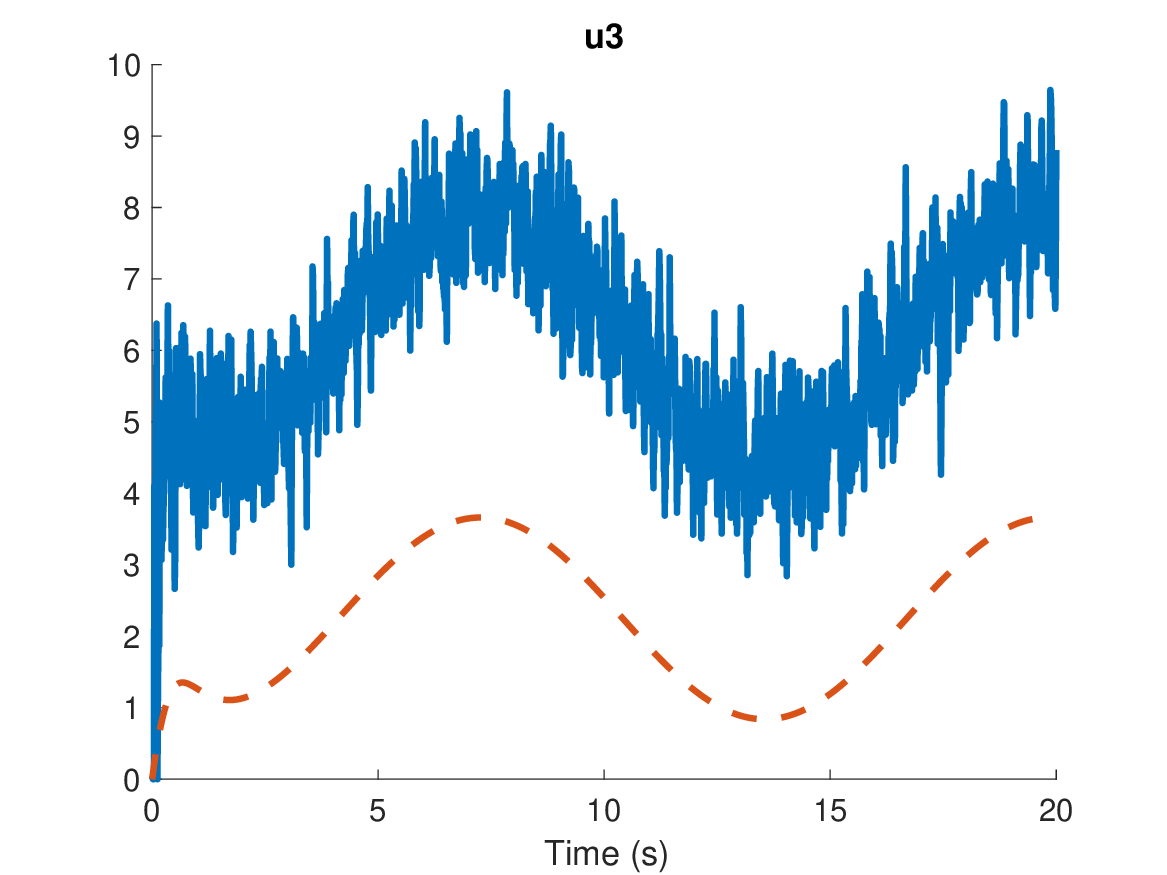,width=0.325\textwidth}}
\caption{Control inputs}\label{Bu}
\end{figure*}

\begin{figure*}[!ht]
\centering
\subfigure[\footnotesize $\theta_1$ (blue--), $\theta_1^\star$ (red - -)]
{\epsfig{figure=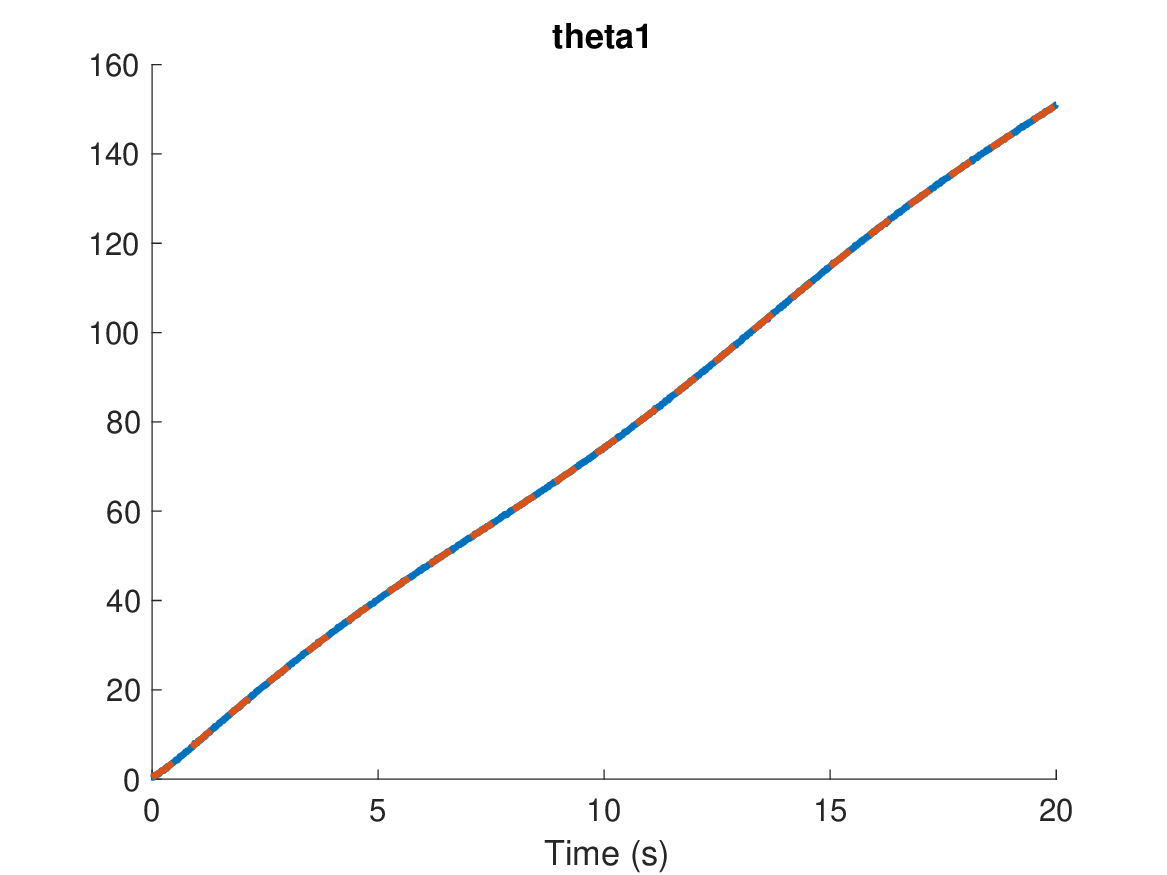,width=0.325\textwidth}}
\subfigure[\footnotesize $\theta_2$ (blue--), $\theta_2^\star$ (red - -)]
{\epsfig{figure=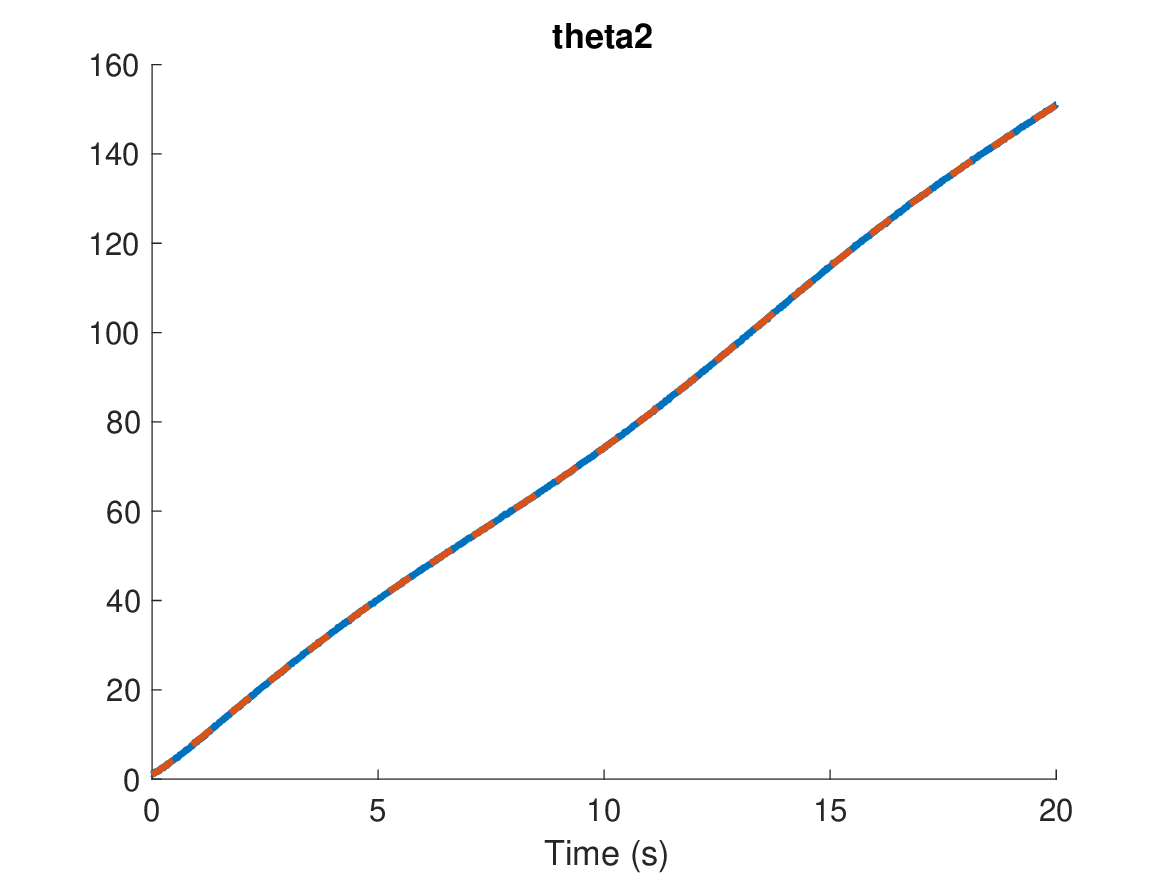,width=0.325\textwidth}}
\subfigure[\footnotesize$\theta_3$ (blue--), $\theta_3^\star$ (red - -)]
{\epsfig{figure=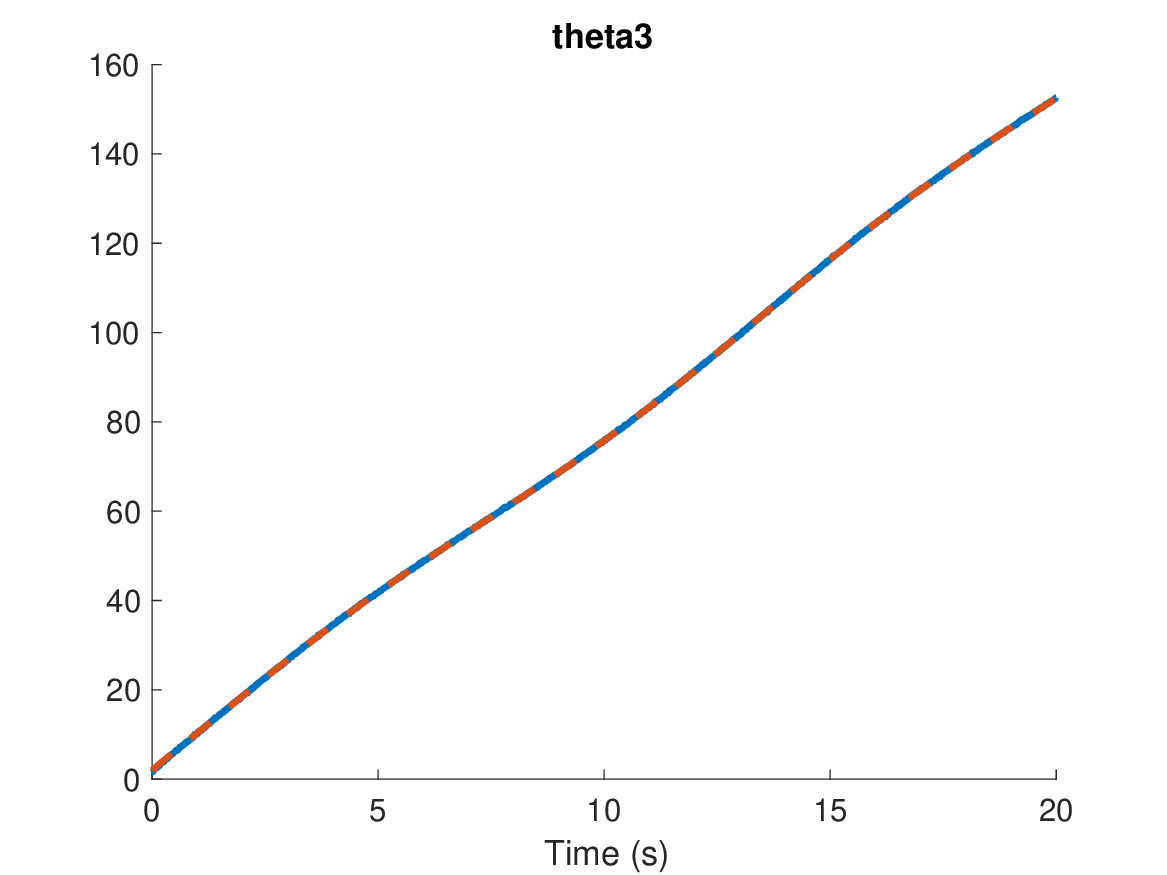,width=0.325\textwidth}}
\caption{Outputs}\label{By}
\end{figure*}

\begin{figure*}[!ht]
\centering
\subfigure[\footnotesize $\dot\theta_1$ (blue--), $\dot\theta_1^\star$ (red - -)]
{\epsfig{figure=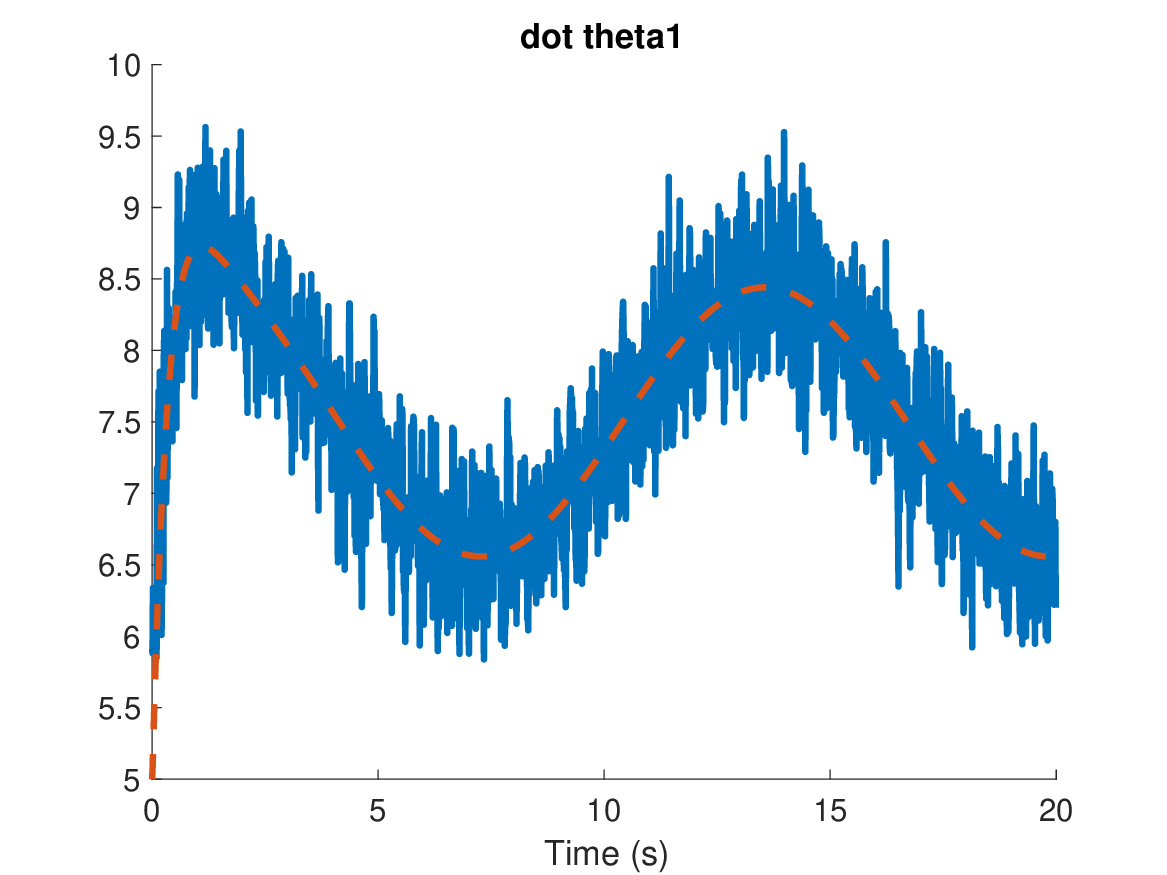,width=0.325\textwidth}}
\subfigure[\footnotesize $\dot\theta_2$ (blue--), $\dot\theta_2^\star$ (red - -)]
{\epsfig{figure=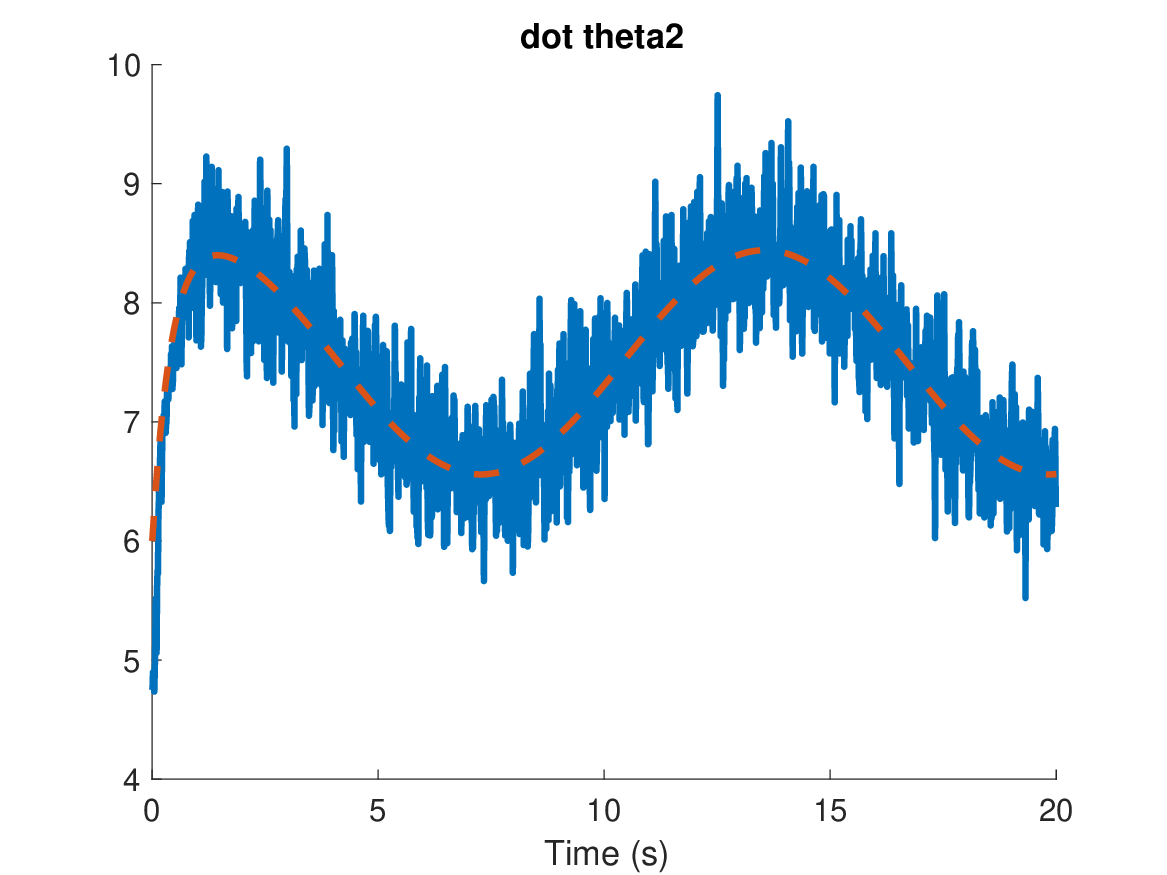,width=0.325\textwidth}}
\subfigure[\footnotesize$\dot\theta_3$ (blue--), $\dot\theta_3^\star$ (red - -)]
{\epsfig{figure=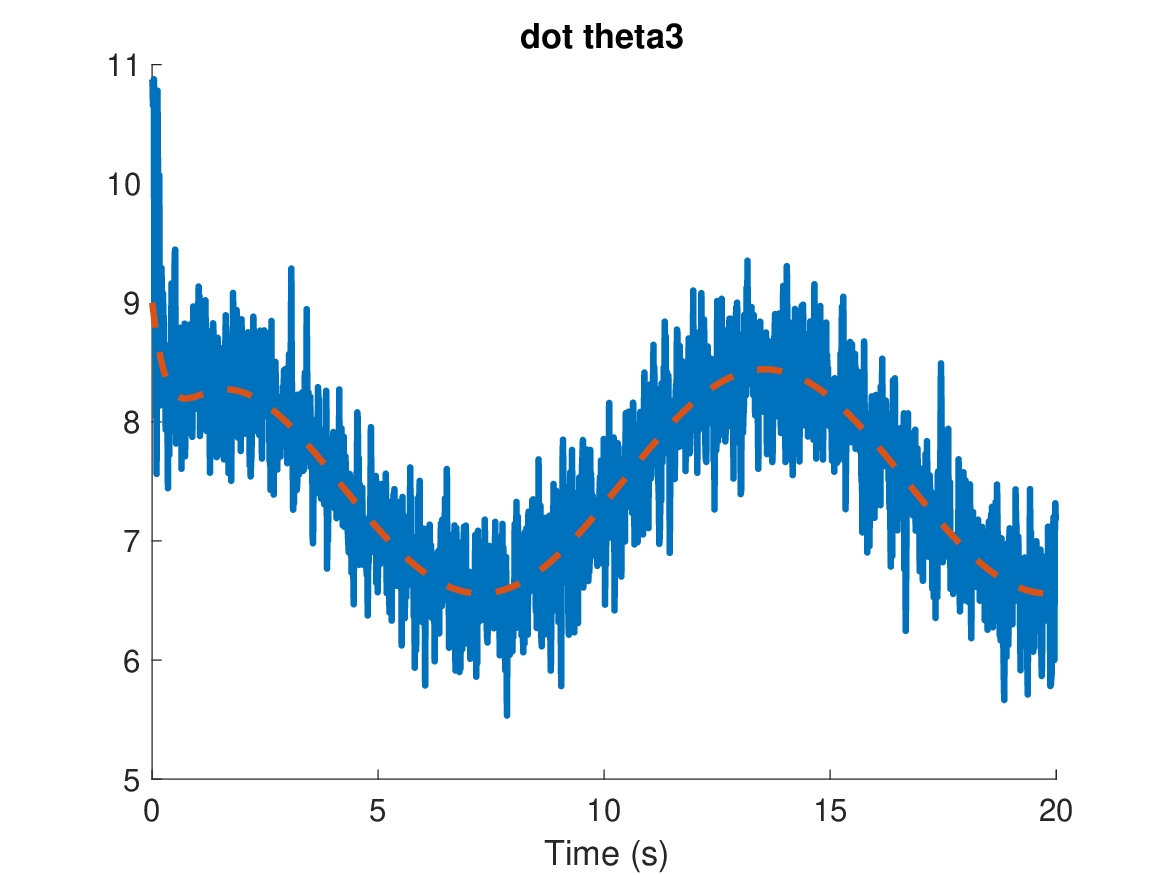,width=0.325\textwidth}}
\caption{Time derivative outputs}\label{Bdy}
\end{figure*}

Consider the case of three oscillators.
Set $N = 3$ in Eq.~\eqref{eq:model:kuramoto3}, and take $\omega_1=5$, $\omega_2=7$, $\omega_3=8$, $K=1$, $K_{P, i}=1$,  $a_{1,2}=a_{1,3}=a_{2,1}=a_{2,3}=a_{3,1}=a_{3,2}=1$. Introduce also the uncertainties $\Delta_1=1.2$, $\Delta_2=0.8$, $\Delta_3=1.2$, $\Delta_4=0.8$, $\Delta_5=0.8$, $\Delta_6=1.2$, $\Delta_7=0.8$. Set $c_1=\pi/2$, $c_2=\pi/2$, $c_3=\pi$, $\theta_1(0)=0.5\Delta_5$, $\theta_2(0)=\Delta_6$, $\theta_3(0)=2\Delta_7$. Figures \ref{By} and \ref{Bdy} depict synchronization towards $f(t) + c_i = 2\sin (0.5 t)+7.5t+7+c_i$:\\
\small
\begin{equation}
\begin{cases}
  \label{eq:model:kuramoto3}
\dot\theta_1 =\omega_1\Delta_1 +u_1 \frac{K\Delta_4}{N}\left(a_{1,2} \sin(\theta_2- \theta_1)+a_{1,3} \sin(\theta_3- \theta_1)\right)\\ 
\dot\theta_2 =\omega_2\Delta_2 +u_2 \frac{K\Delta_4}{N}\left(a_{2,1} \sin(\theta_1- \theta_2)+a_{2,3} \sin(\theta_3- \theta_2)\right) \\
\dot\theta_3 =\omega_3\Delta_3 +u_3 \frac{K\Delta_4}{N}\left(a_{3,1} \sin(\theta_1- \theta_3)+a_{3,2} \sin(\theta_2- \theta_3)\right) 
\end{cases}
\end{equation}\normalsize
The sampling period is $T_e=0.01$\,s. The measures of $\theta_i(t)$, $i = 1, 2, 3$, is corrupted by an additive white Gaussian noise $\mathcal{N}(0,0.1)$. The results may be found in Figures \ref{Bu}, \ref{By}, \ref{Bdy} and \ref{Be}.

\section{Additive control variables}\label{add}
Consider with \cite{Mao16} the following modeling where, contrarily to Eq. \eqref{eq:model:kuramoto}, the control variables are additive:
\begin{equation}
    \label{eq:model:kuramoto:control:ad}
\dot\theta_i=\omega_i +\frac{K}{N}\sum_{j=1}^N a_{i,j} \sin(\theta_j- \theta_i)+u_i
\end{equation}
The above system is again flat: $\theta_1$, $\theta_2$, $\theta_3$ are again flat outputs. Eq. \eqref{eq:model:kuramoto:control:ad} yields
\begin{equation}\label{flat2}
   u_i =  \dot{\theta}_i - \omega_i-\frac{K}{N}\sum_{j=1}^N a_{i,j} \sin(\theta_j - \theta_i)
\end{equation}
Eq. \eqref{flat2} shows that avoiding singularities like in Sect. \ref{sing} becomes pointless. Choosing an appropriate open-loop reference trajectory becomes easy.  
Differentiate Eq.~\eqref{eq:model:kuramoto:control:ad}: 
\begin{equation*}
d\dot\theta_i = 
\frac{K}{N}\sum_{j=1}^N a_{i,j} (d\theta_j - d\theta_i) \cos(\theta_j- \theta_i) + du_i, \quad i = 1,\ldots, N
\end{equation*}
It yields the Eq.~\eqref{eq:homeo} where now
\begin{eqnarray*}
    F_i &=& 
    \frac{K}{N}
    \sum_{j=1}^N a_{i,j}
    (\delta \theta_j - \delta \theta_i)\cos(\theta_{j}^{\star} - \theta_{i}^{\star} )
    \\
    \alpha_i &=& 1, \quad\qquad i = 1,\ldots, N
\end{eqnarray*}

\begin{figure*}[!ht]
\centering
\subfigure[\footnotesize $u_1$ (blue--), $u_1^\star$ (red - -)]
{\epsfig{figure=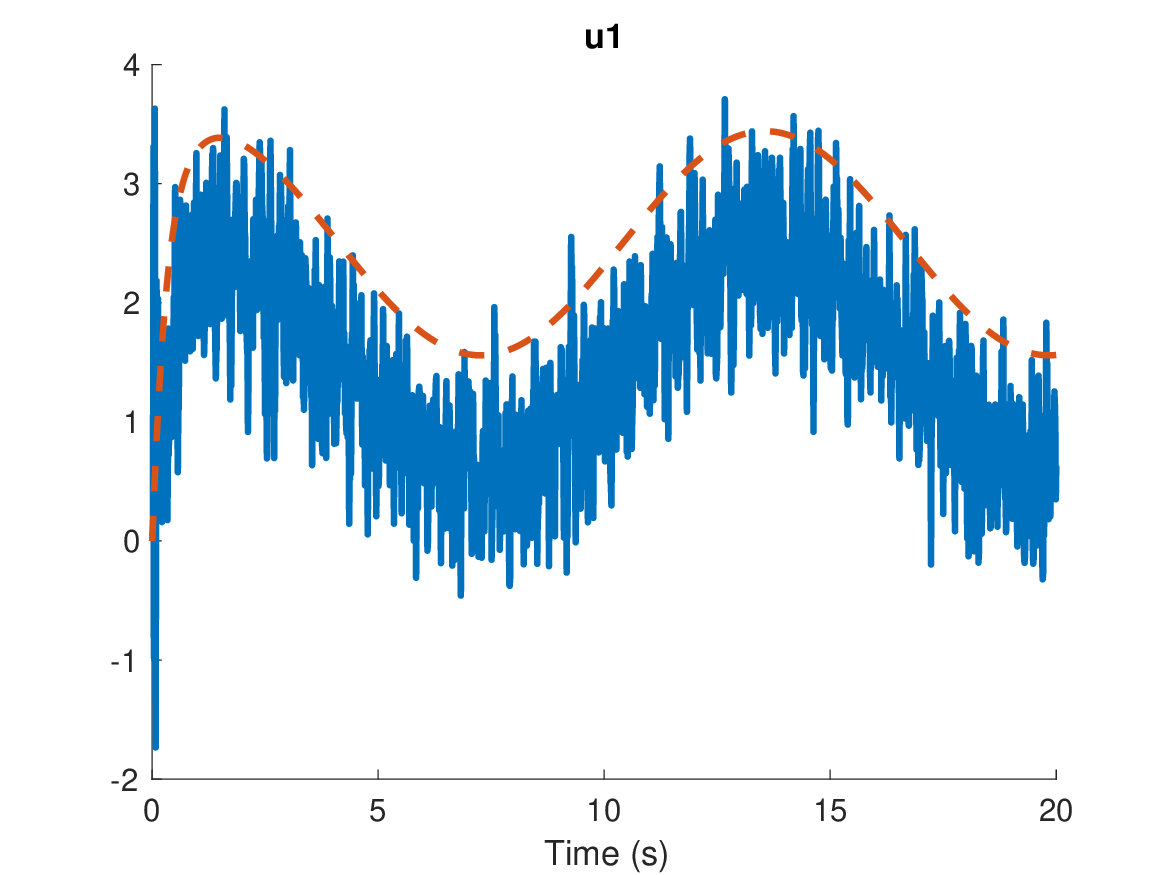,width=0.325\textwidth}}
\subfigure[\footnotesize $u_2$ (blue--), $u_2^\star$ (red - -)]
{\epsfig{figure=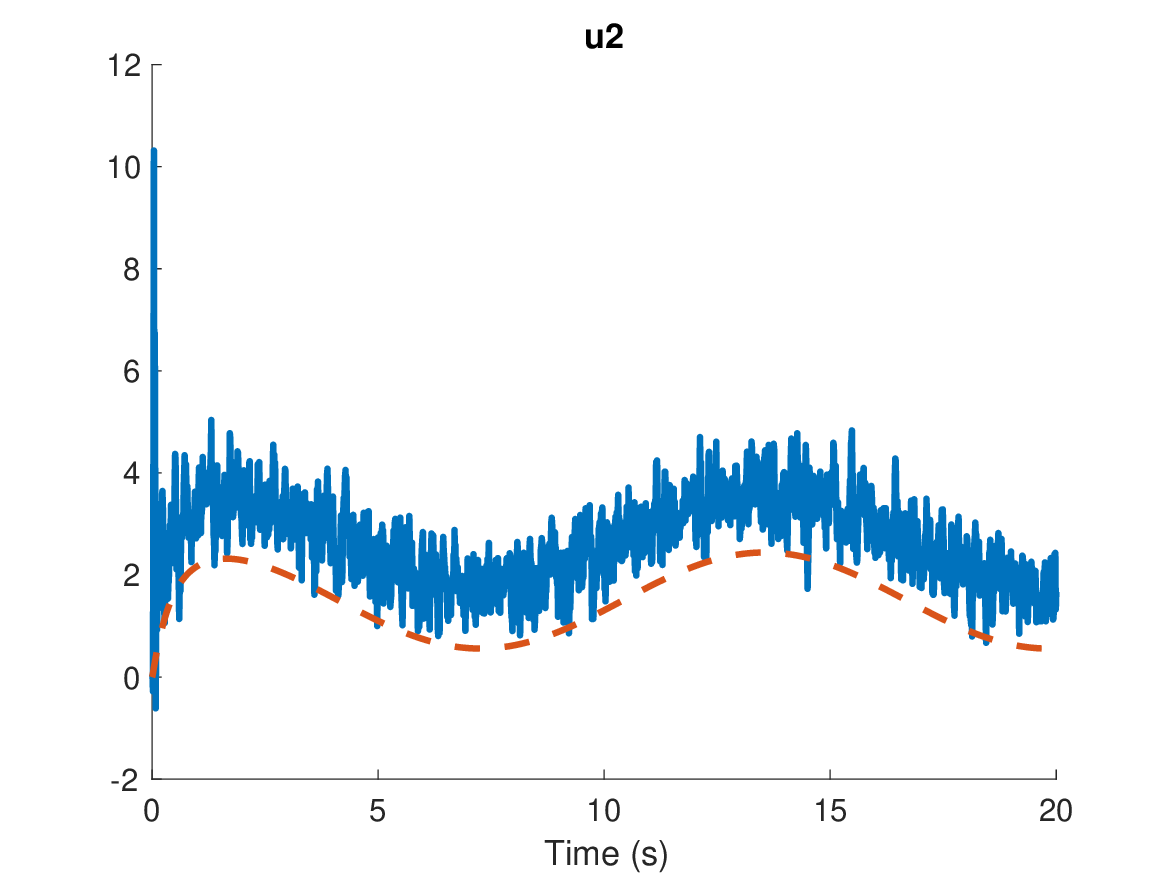,width=0.325\textwidth}}
\subfigure[\footnotesize$u_3$ (blue--), $u_3^\star$ (red - -)]
{\epsfig{figure=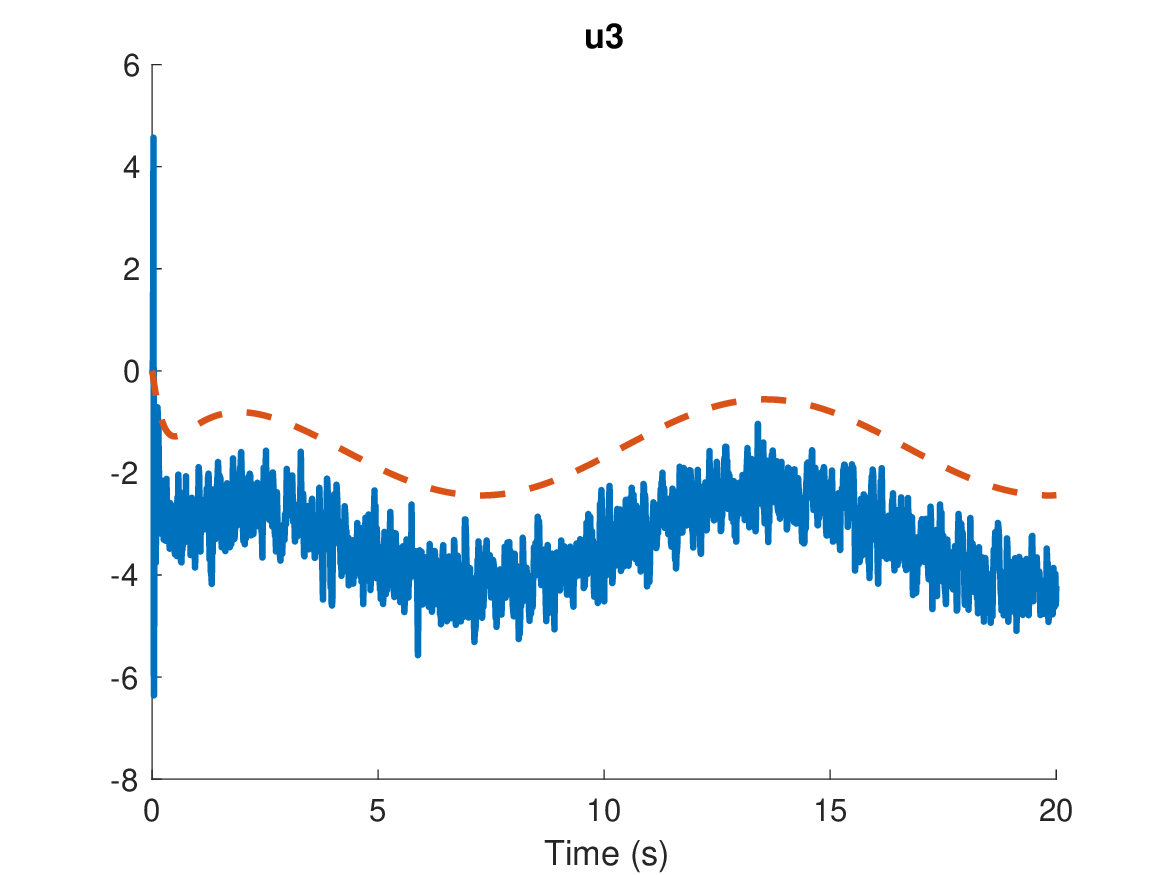,width=0.325\textwidth}}
\caption{Additive case: control inputs}\label{Beu}
\end{figure*}

\begin{figure*}[!ht]
\centering
\subfigure[\footnotesize $\theta_1$ (blue--), $\theta_1^\star$ (red - -)]
{\epsfig{figure=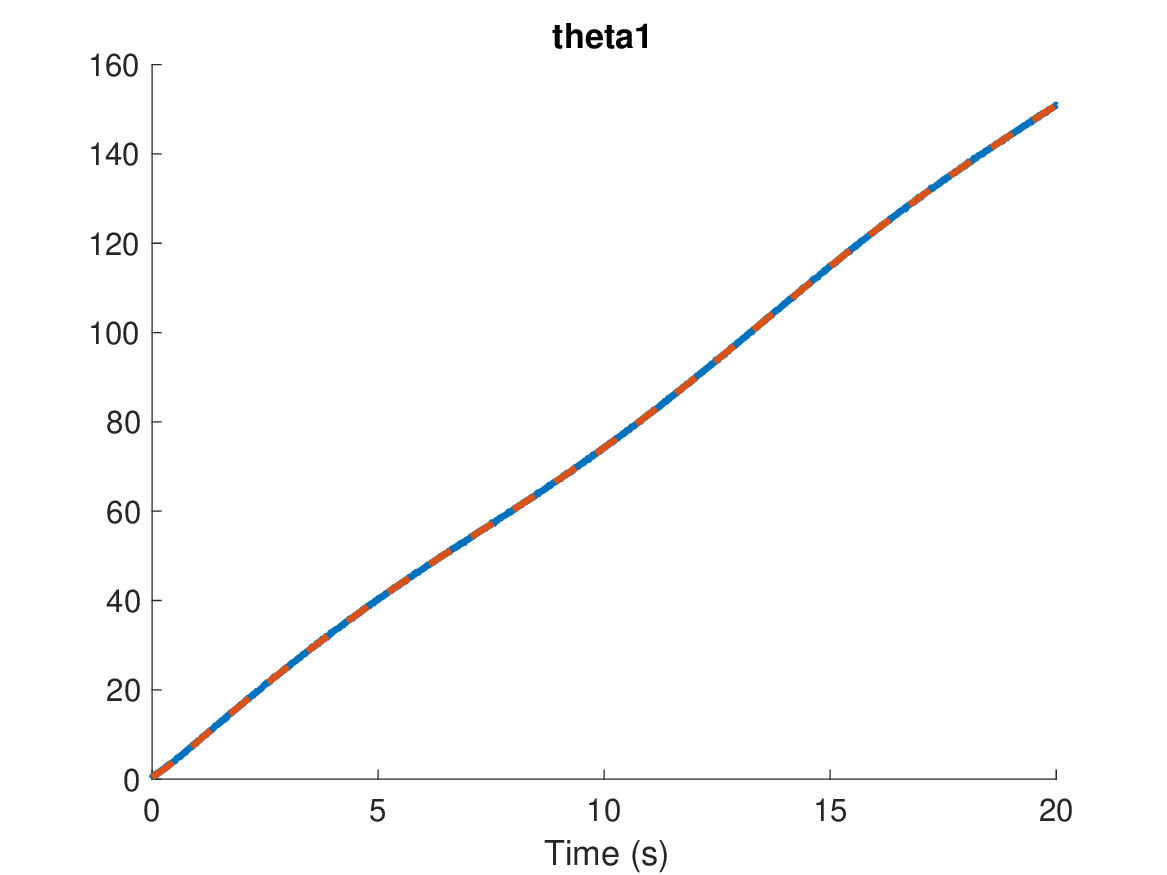,width=0.325\textwidth}}
\subfigure[\footnotesize $\theta_2$ (blue--), $\theta_2^\star$ (red - -)]
{\epsfig{figure=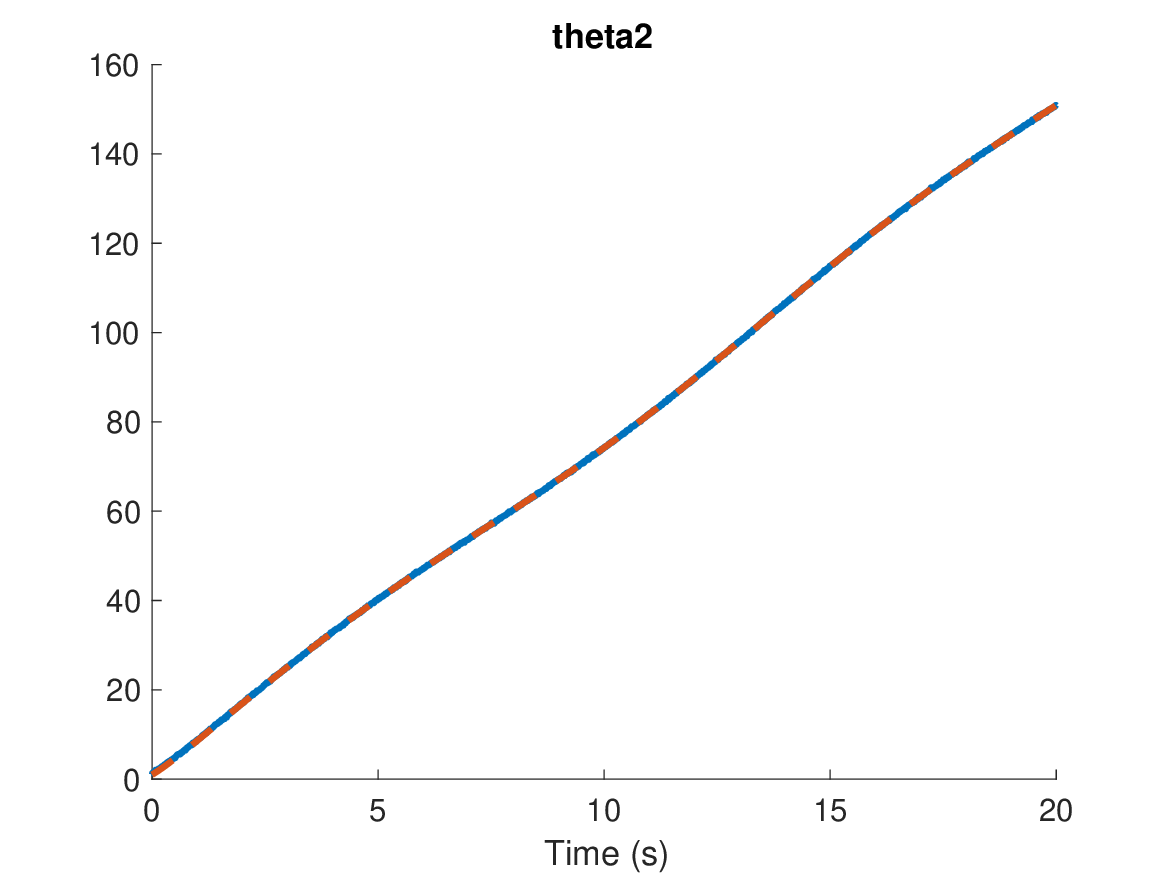,width=0.325\textwidth}}
\subfigure[\footnotesize$\theta_3$ (blue--), $\theta_3^\star$ (red - -)]
{\epsfig{figure=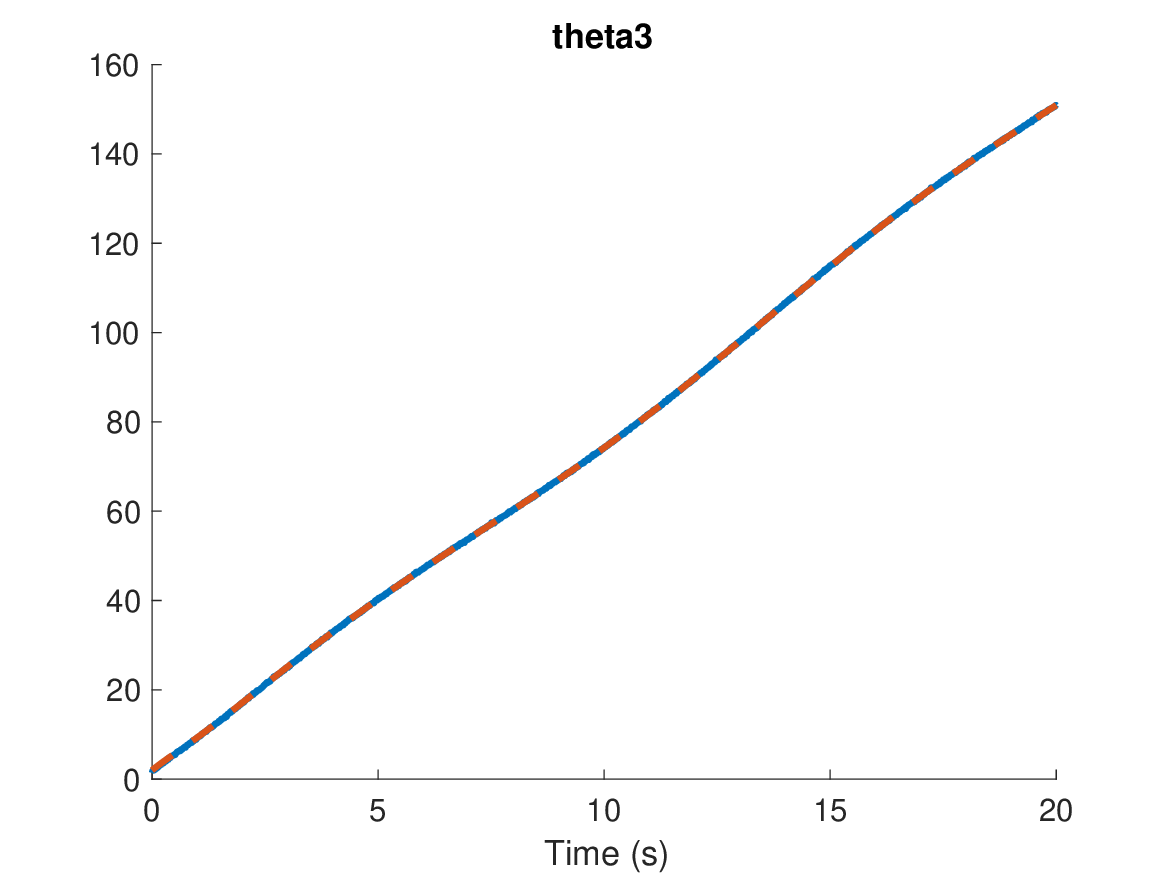,width=0.325\textwidth}}
\caption{Additive case: outputs}\label{Bey}
\end{figure*}

\begin{figure*}[!ht]
\centering
\subfigure[\footnotesize $\dot\theta_1$ (blue--), $\dot\theta_1^\star$ (red - -)]
{\epsfig{figure=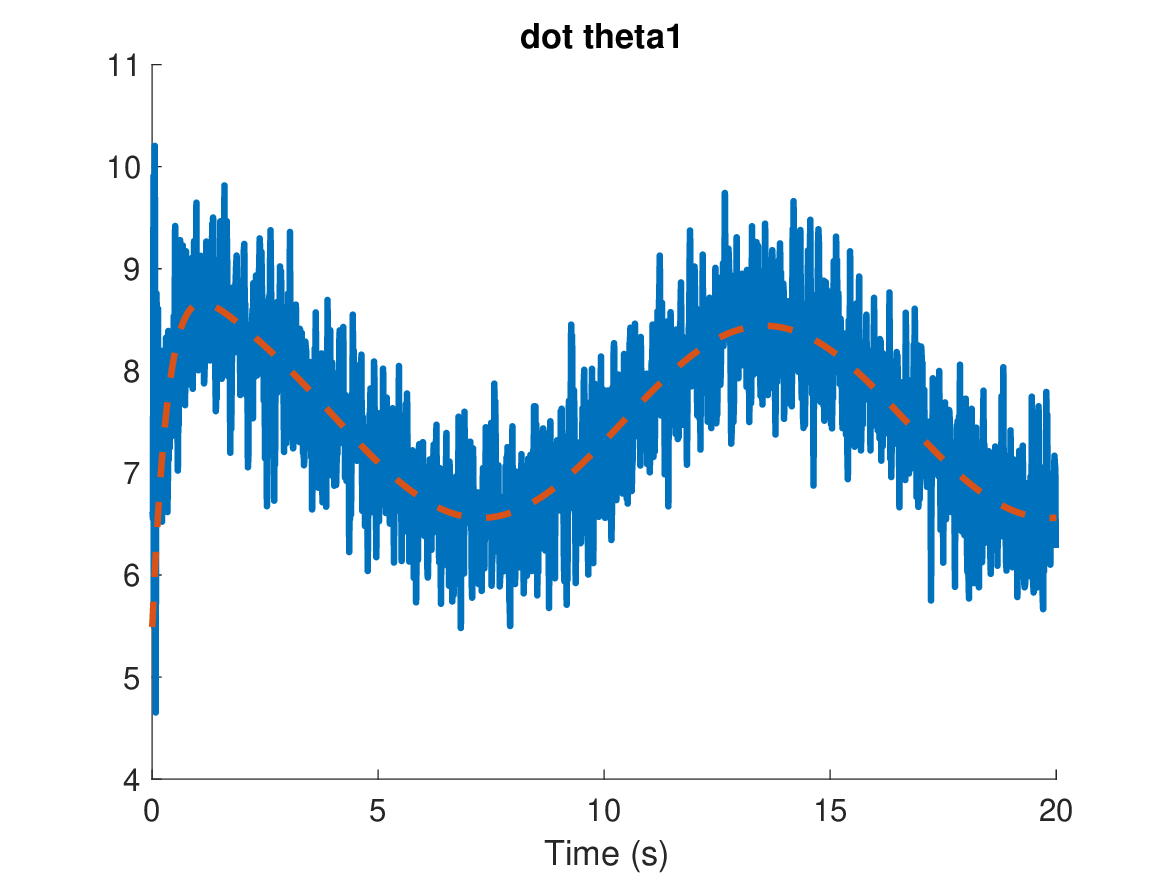,width=0.325\textwidth}}
\subfigure[\footnotesize $\dot\theta_2$ (blue--), $\dot\theta_2^\star$ (red - -)]
{\epsfig{figure=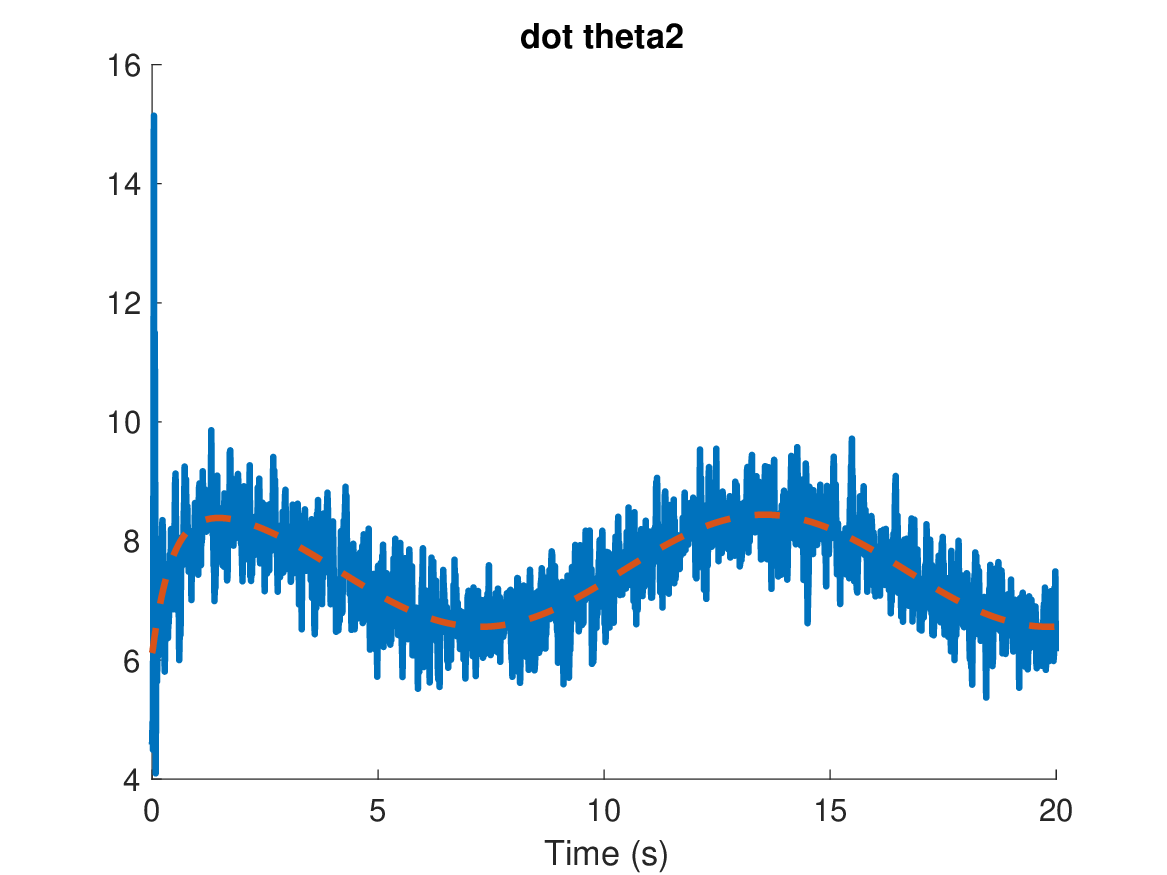,width=0.325\textwidth}}
\subfigure[\footnotesize$\dot\theta_3$ (blue--), $\dot\theta_3^\star$ (red - -)]
{\epsfig{figure=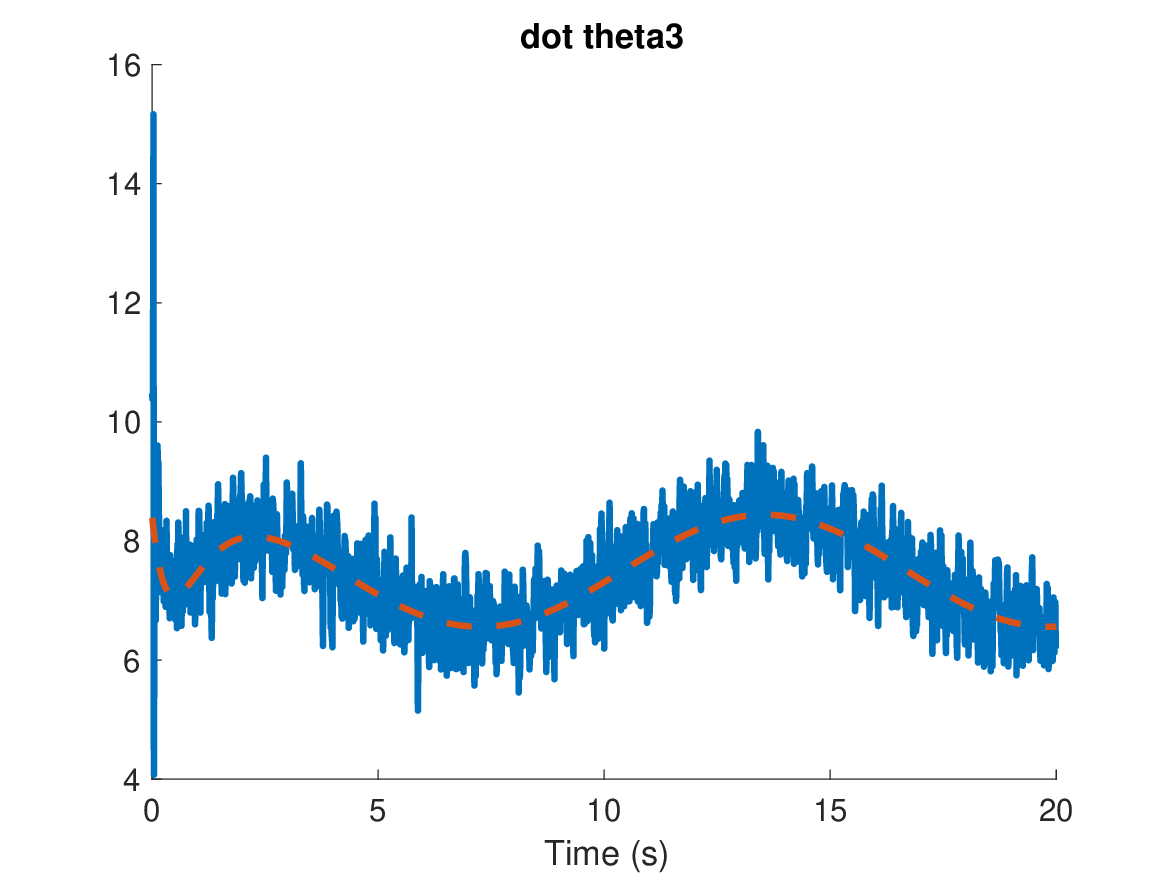,width=0.325\textwidth}}
\caption{Additive case: time derivative outputs}\label{Bedy}
\end{figure*}

\begin{figure*}[!ht]
\centering
\subfigure[\footnotesize $\delta \theta_1$]
{\epsfig{figure=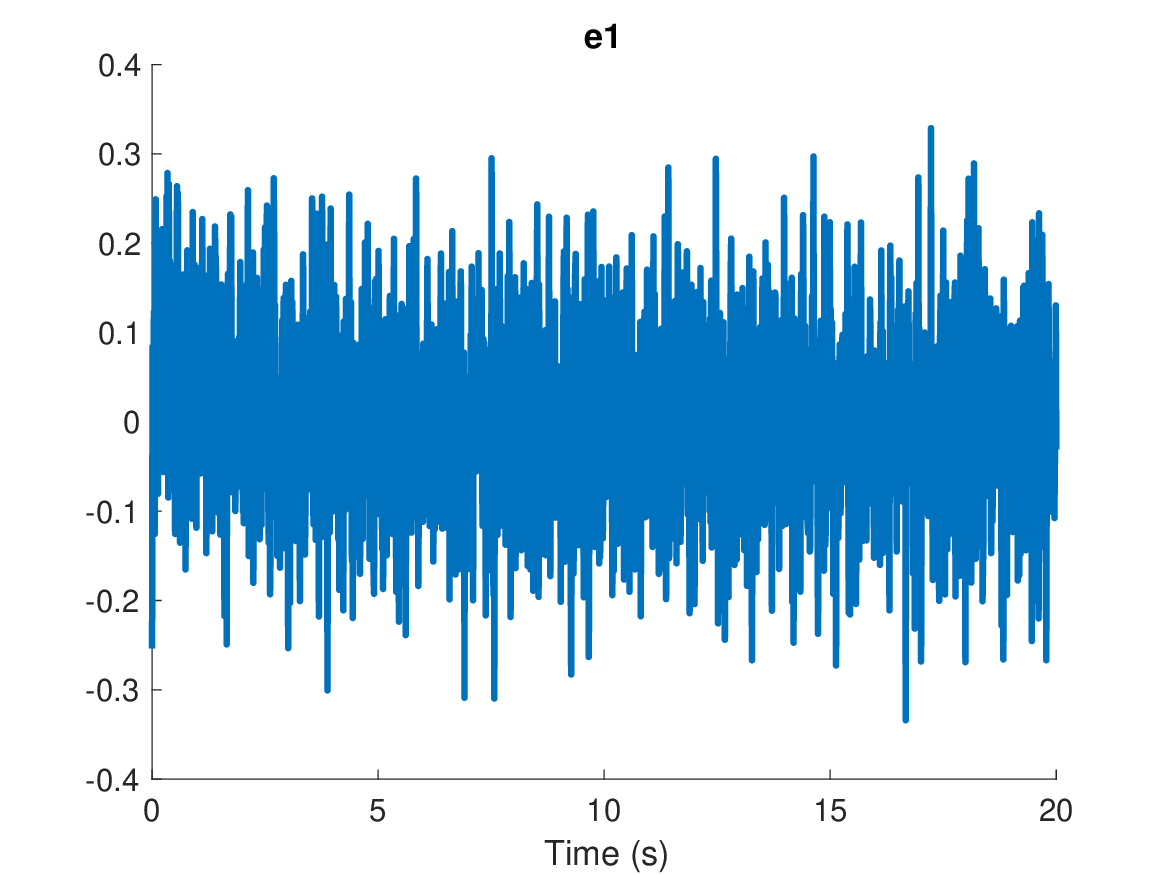,width=0.325\textwidth}}
\subfigure[\footnotesize $\delta \theta_2$]
{\epsfig{figure=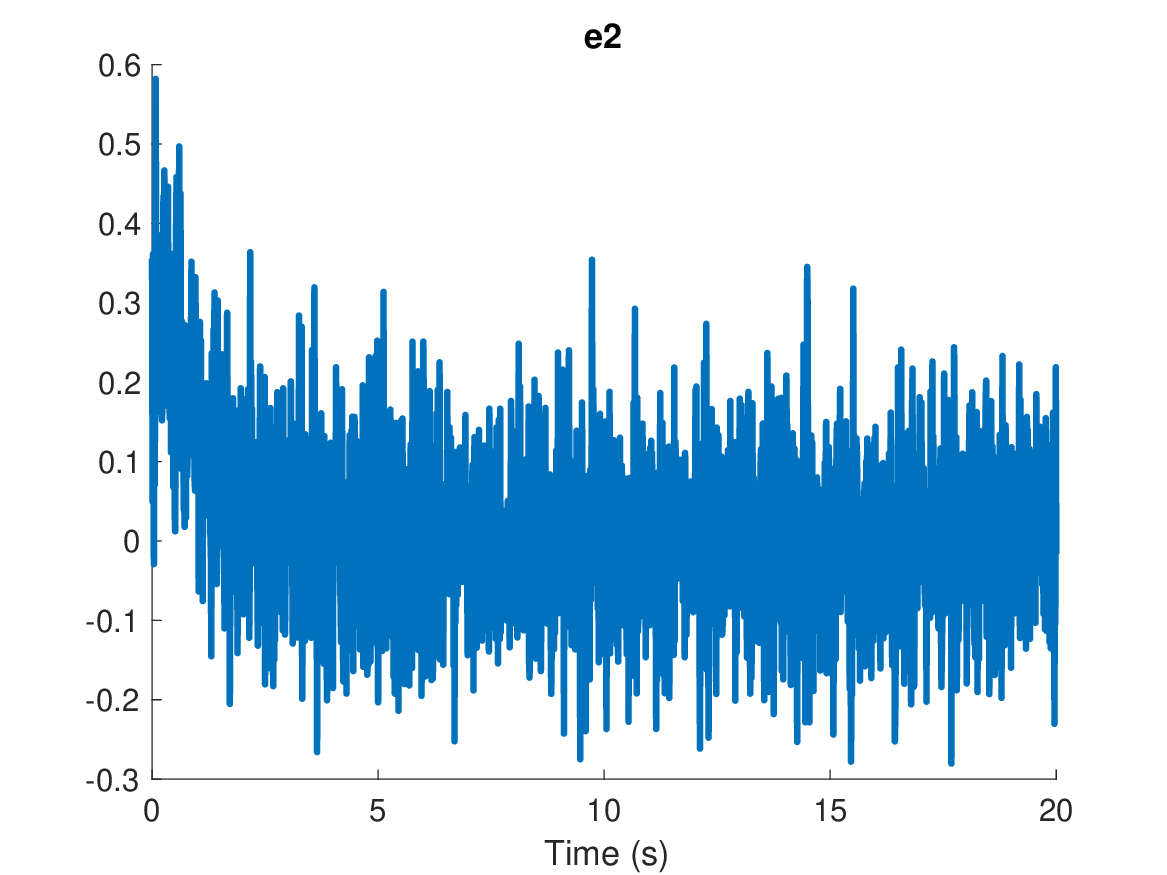,width=0.325\textwidth}}
\subfigure[\footnotesize $\delta \theta_3$]
{\epsfig{figure=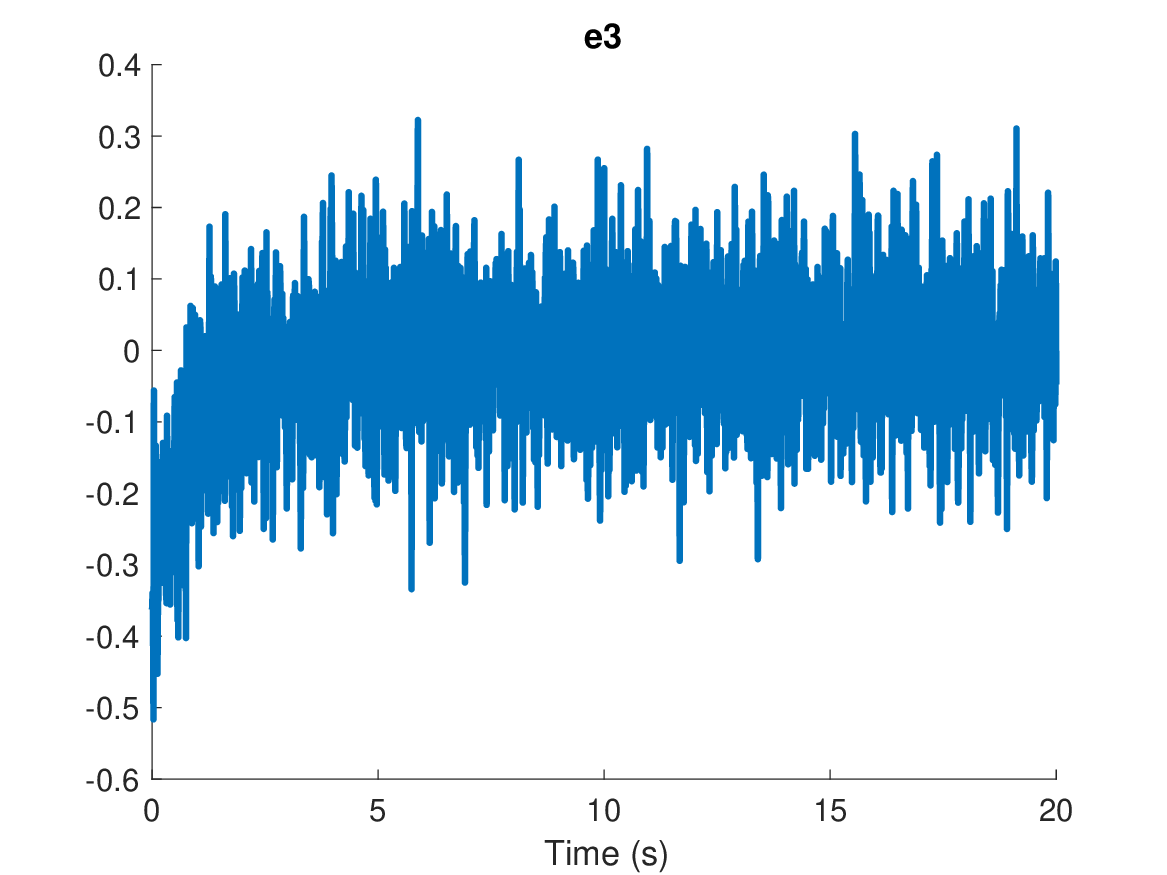,width=0.325\textwidth}}
\caption{Additive case: tracking errors}\label{Bee}
\end{figure*}

With the same uncertainties as in Sect. \ref{exp}, the results obtained with the corresponding homeostat are depicted in Fig. \ref{Beu}, \ref{Bey}, \ref{Bedy}. Fig. \ref{Bee} demonstrates that the tracking error is negligible. Note that we are able to track the same synchronisation function, $c_i=\frac{\pi}{2}$.\\




\clearpage

\section{Conclusion}\label{conclu}
Our approach was here only computer illustrated via a small number of oscillators. To do it with a large number like \cite{Bottcher22} should not be a problem with the type of tools we are developing. But this should of course be confirmed with convincing computer experiments. It is well known (see, e.g., \cite{cohen}, and references therein) that the influence of neurosciences on deep learning is essential. Let us cite here \cite{neurosci} where our methods demonstrate their efficiency in some preliminary investigations about epilepsy. 


Recent advances in control theory, like HEOL, have demonstrated remarkable capabilities in addressing complex system behaviors. While current AI tools show impressive achievements in many areas, the fundamental concepts from control theory could significantly enhance their theoretical foundations and practical applications.\footnote{See \cite{predict} for \emph{model-free predictive control} (\emph{MFPC}), which is related to \emph{reinforcement learning} (\emph{RL}) (see, e.g., \cite{recht}), i.e., a key ingredient of today's AI.} Cybernetics,\footnote{See \cite{wiener}. The first edition was published in 1948.} which ought to be considered as one of the historical roots of~AI, highlights the importance of feedback loops, a concept that remains challenging to implement via ANNs (\cite{herzog}). In this perspective, we propose viewing control theory as an intrinsic component of AI, following the natural legacy of Wiener's cybernetics. See, e.g., \cite{kline} and \cite{najim} for a glimpse of the fascinating historical background.
%





\begin{thebibliography}{3}
\providecommand{\natexlab}[1]{#1}
\providecommand{\url}[1]{\texttt{#1}}
\providecommand{\urlprefix}{URL }
\expandafter\ifx\csname urlstyle\endcsname\relax
  \providecommand{\doi}[1]{doi:\discretionary{}{}{}#1}\else
  \providecommand{\doi}{doi:\discretionary{}{}{}\begingroup
  \urlstyle{rm}\Url}\fi


\bibitem[Acebr\'{o}n et al.(2005)]{Acebron05}
Acebr\'{o}n, J. A., Bonilla, L. L., P\'{e}rez Vicente, C. J., Ritort, F., Spigler, R. (2005). The Kuramoto model: A simple paradigm for synchronization phenomena. {\it Rev. Mod. Phys.}, 77, 137-185.

\bibitem[Ashby(1960)]{ashby}
W.R. Ashby (1960). {\it Design for a Brain}. Chapman \& Hall.

\bibitem[Bensoussan et al.(2022)]{ben}
Bensoussan, A., Li, Y., Nguyen, D.P.C., Tran, M.-B., Chi, S., Yam, P., Zhou, X. (2022). Machine learning and control theory.
Trélat, E., Zuazua, E. (Eds):
{\it Handbook of Numerical Analysis}, Elsevier, vol. 23, pp. 531-558.

\bibitem[{B\"ottcher et~al.(2022)Böttcher, Antulov-Fantulin, and
  Asikis}]{Bottcher22}
Böttcher, L., Antulov-Fantulin, N., Asikis, T. (2022).
\newblock {AI} {P}ontryagin or how artificial neural networks learn to control
  dynamical systems.
\newblock \textit{Nature Commun.}, 13, 333.


\bibitem[Breakspear et al.(2010)]{Breakspear10}
Breakspear, M., Heitmann, S. Daffertshofer, A. (2010).
\newblock Generative models of cortical oscillations: neurobiological implications of the Kuramoto model.
\newblock \textit{Front. Human Neurosci.}, 4, 190.

\bibitem[Cerf and Rutten(2023)]{cerf}
Cerf, S., Rutten, E. (2023).
Combining neural networks and control: potentialities, patterns and perspectives. {\it IFAC-PapersOnLine},
56, 9036-9049.

\bibitem[Chopra and Spong(2005)]{chopra}
Chopra, N., Spong, M.W. (2005). On Synchronization of Kuramoto Oscillators. {\it 44th IEEE Conf. Decis. Contr.}. Seville, pp. 3916-3922. 

\bibitem[Cohen et al.(2022)]{cohen}
Cohen, Y., Engel, T.A., Langdon, C., Lindsay, G.W., Ott, T.,
Peters, M.A.K.,   Shine, J.M., Breton-Provencher, V., Ramaswamy, S. (2022). 
Recent advances at the interface of neuroscience and artificial neural networks. {\it J. Neurosci.},  9, 8514-8523.


\bibitem[Delaleau and Hagenmeyer(2002)]{DelalHagen02jesa}
Delaleau, E., Hagenmeyer, V. (2002).
Commande prédictive non linéaire fondée sur la 
platitude du moteur à induction: 
Application au positionnement de précision.
{\it J. Europ. Syst. Automat.}, 36, 737-748.

\bibitem[Dev et al. (2021)]{dev}
Dev, P., Jain, S., Arora, P.K., Harish Kumar, H. (2021).
Machine learning and its impact on control systems: A review.
{\it Materials Today Proc.}, 47, 3744-3749.

\bibitem[D\"{o}rfler and Bullo (2014)]{Dorfler14}
D\"{o}rfler, F., Bullo, F. (2014). Synchronization in complex networks of phase oscillators: A survey. \textit{Automatica}, 50, 1539-1564.


\bibitem[D\"{o}rfler et al. (2013)]{Dorfler13}
D\"{o}rfler, F., Chertkov, M., Bullo, F. (2014). Synchronization in complex oscillator networks and smart grids. {\it Proc. Natl. Acad. Sci. U S A.}, 110, 2005-2010. 

\bibitem[Fliess and Join(2013)]{mfc1}
Fliess, M., Join, C. (2013). Model-free control. {\it Int. J. Contr.}, 86, 2228-2252.

\bibitem[Fliess and Join(2021)]{vancouver}
Fliess, M., Join, C. (2021). Machine learning and control engineering: The model-free case. Arai, K., Kapoor, S., Bhatia, R. (eds): {\it Proc. Future Techno. Conf. (FTC) 2020}, Advances in Intelligent Systems and Computing, vol. 1288, pp. 258-278. Springer.

\bibitem[Fliess and Join(2022)]{mfc2}
Fliess, M., Join, C. (2022). An alternative to proportional-integral and proportional-integral-derivative regulators: Intelligent proportional-derivative regulators. {\it Int. J. Robust Nonlin. Contr.}, 32, 9512-9524.

\bibitem[Fliess et al. (1995)]{flmr_ijc}
Fliess, M., L{\'e}vine, J., Martin, P., Rouchon, P. (1995). Flatness and defect of non-linear systems: introductory theory and examples, {\it Int. J. Contr.}, 61, 1327-1361.

\bibitem[Fliess et al.(1999)]{flmr_ieee}
Fliess, M., L\'{e}vine, J., Martin, P., Rouchon, P. (1999).
A Lie-B\"acklund approach to equivalence and flatness
of nonlinear systems, {\it IEEE Trans. Automat. Contr.}, 44, 922-937.

\bibitem[Franci et al.(2012)]{lss}
Franci, A., Chaillet, A., Panteley, E., Lamnabhi-Lagarrigue, F. (2012). Desynchronization and inhibition of Kuramoto oscillators by scalar mean-field feedback. {\it Math. Contr. Sign. Syst.}, 24, 167-217.


\bibitem[{Ha et al.(2016a)}]{Ha16}
Ha, S.Y., H.K. Kim, S.W. Ryoo (2016a).
\newblock Emergence of phase-locked states for the Kuramoto
model in a large coupling regime.
\newblock \emph{Commun. Math. Sci.}, 14(4) p.~1073–-1091.
\newblock \doi{https://doi.org/10.1038/s41467-021-27590-0}.

\bibitem[Ha et al.(2016b)]{HaSIAM}
Ha, S.Y., Noh, S.E., Park, J. (2016b). Synchronization of Kuramoto oscillators with adaptive couplings. {\it SIAM J. App. Dyn. Syst.}, 15, 162-194.

\bibitem[Hagenmeyer(2003)]{hagen}
Hagenmeyer (2003). {\it Robust Nonlinear Tracking Control Based on Differential Flatness}. VDI Verlag.

\bibitem[Herzog et al.(2020)]{herzog}
Herzog, S., Tetzlaff, C., W\"{o}rg\"{o}tter, F. (2020). Evolving artificial neural networks with feedback.
{\it Neural Netw.}, 123, 153-162.

\bibitem[Jin et al.(2020)]{jin}
Jin, W., Wang, Z., Yang, Z., Mou, S. (2020). Pontryagin differentiable programming:
An end-to-end learning and control framework. {\it 34th Conf. Neural Informat. Process. Syst. (NeurIPS)}, Vancouver. 

\bibitem[Join et al.(2013)]{chaxel}
Join, C., Chaxel, F., Fliess, M. (2013). ``Intelligent'' controllers on cheap and small programmable devices. {\it 2nd Int. Conf. Contr. Fault-Tolerant Syst.} , Nice.

\bibitem[Join et al.(2024a)]{heol}
Join, C., Delaleau, E., Fliess, M. (2024a). Flatness-based control revisited: The {\it HEOL} setting, \textit{C.R. Math.}, 362, 1693-1706.

\bibitem[Join et al.(2025)]{predict}
Join, C., Delaleau, E., Fliess, M. (2025).
\newblock Model-free predictive control: Introductory algebraic calculations, and a brief comparison with HEOL.
\textit{Submitted}.


\bibitem[Join et al.(2024b)]{neurosci}
Join, C., Jovellar, D.B., Delaleau, E., Fliess, M. (2024b).
Detection and suppression of epileptiform seizures via model-free control and derivatives in a noisy environment. {\it 12th Int. Conf. Syst. Contr.}, Batna.

\bibitem[Kline(2011)]{kline}
Kline, C. (2011). Cybernetics, automata studies, and the Dartmouth conference on artificial intelligence, {\it IEEE Ann. History Comput.}, 33, 5-16.

\bibitem[Kuramoto (1975)]{Kuramoto75}
Kuramoto, Y. (1975).  Self-entrainment of a population of coupled non-linear oscillators. Araki, H. (Ed.) {\it International Symposium on Mathematical Problems in
Theoretical Physics}, pp. 420–422, Springer.

\bibitem[Kuramoto (1984)]{Kuramoto84}
Kuramoto, Y. (1984) {\it Chemical Oscillations, Waves and Turbulence}. Springer.

\bibitem[LeCun et al.(2015)]{cun}
LeCun, Y., Bengio, Y., Hinton, G. (2015).
Deep learning, \textit{Nature}, 521, 436-444.

\bibitem[L\'{e}vine (2009)]{levine}
L\'{e}vine, J. (2009). \textit{Analysis and Control of Nonlinear Systems: A Flatness-based Approach}. Springer.

\bibitem[{Mao and Zhang(2016)}]{Mao16}
Mao, Y., Zhang, Z.(2016).
\newblock Distributed frequency synchronization and phase-difference tracking
for Kuramoto oscillators and its application to islanded microgrids.
\newblock \emph{IEEE 55th Conf. Decis. Contr.}. Las Vegas. pp. 4364-4369. 
\bibitem[Miller et al.(1991)]{mit}
Miller, W.T., Werbos, P.J., Sutton, R.S. (eds) (1991). \textit{Neural Networks for Control}. MIT Press.

\bibitem[{Najim(2024)}]{najim}
Najim, M. (2024), Editorial: Did Norbert Wiener’s cybernetics crystallize in France and incubate AI in China? {\it Kybernetes}, 53, 6155-6165.

\bibitem[{Narendra and Parthasarathy(1990)}]{nar}
Narendra, K.S., Parthasarathy, K. (1990). Identification and control of dynamical systems using neural networks. {\it IEEE Trans. Neural Netw.}, 1, 4-27. 

\bibitem[Recht(2019)]{recht}
Recht, B. (2019).
A tour of reinforcement learning: The view from continuous control.
\textit{Annual Rev. Contr. Robot. Auton. Syst.},
2, 253-279.


\bibitem[Rudolph(2021)]{rudolph}
Rudolph, J. (2021). {\it Flatness-Based Control: An Introduction}, Shaker Verlag.

\bibitem[Sarangapani(2018)]{sara}
Sarangapani, J. (2018). {\it Neural Network Control of Nonlinear Discrete-Time Systems}. CRC Press.

\bibitem[Sira-Ram\'{\i}rez and Agrawal(2004)]{sira}
Sira-Ram\'{\i}rez, H., Agrawal, S.K. (2004). \textit{Differentially Flat Systems}. Marcel Dekker.

\bibitem[Suykens et al.(2010)]{suy}
Suykens, J.A.K., Vandewalle, J.P.L., Bart L.R. De Moor, B.L.R. (2010). {\it Artificial Neural Networks for Modelling and Control of Non-Linear Systems}. Kluwer.

\bibitem[Sutton(1988)]{sutton}
Sutton, R.S. (1988). Artificial intelligence as a control
problem: Comments on the relationship between machine learning and intelligent control. {\it Proc.
IEEE Int. Symp. Intel. Contr.}, pp. 500–507, Arlington.

\bibitem[Strogatz(2003)]{Stroglatz05}
Strogatz, S. (2003). {\it Sync: The emerging science of spontaneous order}. Hyperion.

\bibitem[Wiener(2019)]{wiener}
Wiener, N. (2019). \textit{Cybernetics}. MIT Press (1st ed.: 1948).

\bibitem[Zhou et al.(2024)]{zhou}
Zhou, K., Mao, B., Zhang, Y., Chen, Y., Xiang, Y., Yu, Z., Hao, H., Tang, W., Li, Y., Liu, H., Wang, X., Wang, X. and Qu, J. (2024), A Cable-Actuated Soft Manipulator for Dexterous Grasping Based on Deep Reinforcement Learning. {\it Adv. Intell. Syst.}, 2400112 



\end{thebibliography}
\end{document}